\documentclass[review,onefignum,onetabnum]{siamart190516}
\nolinenumbers

\usepackage{tikz,pgfplots,datatool,enumitem}
\usepackage[tikz]{bclogo}
\usetikzlibrary{plotmarks}
\usepackage{amsmath,amssymb}
\usepackage{color,graphicx,fancybox,pause,fancyvrb}
\usepackage{mathrsfs,arydshln}

\newtheorem{rema}[theorem]{Remark}
\newcommand{\theo}[1]{\begin{theorem}#1\end{theorem}}
\newcommand{\defi}[1]{\begin{definition}#1\end{definition}}
\newcommand{\coro}[1]{\begin{corollary}#1\end{corollary}}
\newcommand{\lemm}[1]{\begin{lemma}#1\end{lemma}}
\newcommand{\remark}[1]{\begin{rema}\rm #1\end{rema}}
\newcommand{\arc}{\arraycolsep0.3ex}

\newcommand{\skipthis}[1]{}

\newcommand{\enu}[1]{\begin{enumerate}[leftmargin=3ex,topsep=0ex] #1 \end{enumerate}}

\newcommand{\mun}[1]{\begin{multline*}#1\end{multline*}}

\newcommand{\eqn}[1]{$$#1$$}
\newcommand{\equ}[1]{\begin{equation}#1\end{equation}}

\newcommand{\gat}[1]{\begin{gather}#1\end{gather}}

\newcommand{\arr}[2]{\begin{array}{#1}#2\end{array}}

\newcommand{\mas}[2]{\left[\begin{array}{#1}#2\end{array}\right]}
\newcommand{\mat}[2]{\left(\begin{array}{#1}#2\end{array}\right)}

\newcommand{\cen}[1]{\begin{center}#1\end{center}}
\newcommand{\epro}{ \mbox{\rule{2mm}{2mm}} }

\newcommand{\red}[1]{{\color{sred}#1}}

\definecolor{yel}{rgb}{1,.9,.6}
\definecolor{lmag}{rgb}{0,0.82,1}

\definecolor{mag}{rgb}{0,0.7,.9}
\definecolor{magenta}{rgb}{1,0,1}
\definecolor{lgrey}{rgb}{0.85,0.85,0.85}
\definecolor{lgr}{rgb}{.9,0.9,0.9}
\definecolor{lgreen}{rgb}{0,0.92,0.7}
\colorlet{lred}{red!20!white}
\definecolor{lblue}{rgb}{1,0.95,0.7}

\definecolor{sblue}{rgb}{.1,.2,.5}
\colorlet{lsblue}{sblue!20!white}
\definecolor{sred}{rgb}{.85,.25,0}
\colorlet{lsred}{sred!20!white}
\definecolor{lye}{rgb}{1,0.95,0.7}

\usepgflibrary{fpu}
\usetikzlibrary{positioning,calc,arrows,shapes,shadows,matrix,backgrounds,external}

\tikzset{
auto,
poi/.style={
minimum size=0,
inner sep=0
},
sys/.style 2 args={
rectangle,
draw,
rounded corners,
drop shadow,
draw=black,
top color=black!20,bottom color=black!0,
minimum height=#2,
minimum width=#1,
inner sep=\dn},
syse/.style 2 args={
rectangle,
draw=none,
minimum height=#2,
minimum width=#1,
inner sep=\dn},
nod/.style={
circle,
draw,
fill=white,
minimum size=5ex
},
sum/.style={circle,draw,draw=black,inner sep=0mm,minimum size=2mm,drop shadow,fill=white,
draw=black!100,top color=black!20,bottom color=black!0},
jun/.style={circle,draw,draw=black,inner sep=0mm,minimum size=0mm},
>={latex},
every path/.style={rounded corners},
lin/.style={color=black,draw,->},
}

\newcommand{\tio}[4]{\coordinate (#1) at ($(#2.south #3)!#4!(#2.north #3)$)}

\def\dn{1ex}
\def\dl{3*\dn}

\tikzstyle{sy0}=[sys={0*\dn}{0*\dn}]
\tikzstyle{sy1}=[sys={12*\dn}{8*\dn}]
\tikzstyle{sy2}=[sys={8*\dn}{6*\dn}]
\tikzstyle{sy3}=[sys={5*\dn}{5*\dn}]
\tikzstyle{sy0}=[sys={0*\dn}{0*\dn}]
\tikzstyle{sye1}=[syse={12*\dn}{8*\dn}]
\tikzstyle{sye2}=[syse={8*\dn}{6*\dn}]
\tikzstyle{sye3}=[syse={5*\dn}{5*\dn}]

\tikzstyle{sy4}=[sys={12*\dn}{6*\dn}]
\tikzstyle{sye4}=[syse={12*\dn}{6*\dn}]

\newcommand{\cge}{\succcurlyeq}
\newcommand{\cle}{\preccurlyeq}
\newcommand{\cl}{\prec}
\newcommand{\cg}{\succ}
\newcommand{\rk}{\text{rank}}
\newcommand{\col}{\text{col}}

\renewcommand{\c}[1]{ {\cal #1} }
\newcommand{\bm}[1]{\text{\boldmath$ #1 $} }

\newcommand{\ab}[1]{\langle #1 \rangle}

\newcommand{\argmin}{ \operatornamewithlimits{arg\,min} }
\newcommand{\diag}{ \operatornamewithlimits{diag} }

\renewcommand{\t}[1]{ \tilde{#1} }

\newcommand{\la}{\lambda}
\newcommand{\La}{\Lambda}
\newcommand{\ga}{\gamma}
\newcommand{\eps}{\varepsilon}

\newcommand{\C}{ {\mathbb{C}} }

\newcommand{\Q}{ {\mathbb{Q}} }
\newcommand{\R}{ {\mathbb{R}} }

\renewcommand{\S}{ {\mathbb{S}} }

\newcommand{\N}{ {\mathbb{N}} }
\newcommand{\Z}{ {\mathbb{Z}} }

\newcommand{\T}{\c{T}}

\renewcommand{\ker}{\text{ker}}

\renewcommand{\P}{ {\mathbb{P}} }
\newcommand{\st}{\ |\ }

\newcommand{\bul}{\bullet}

\newcommand{\te}[1]{\text{\ \ #1\ \ }}

\newcommand{\al}{\alpha}
\newcommand{\be}{\beta}
\newcommand{\si}{\sigma}

\newcommand{\Ac}{\green{A_c}}
\newcommand{\Bc}{\green{B_c}}
\newcommand{\Cc}{\green{C_c}}
\newcommand{\Dc}{\green{D_c}}

\newcommand{\Aq}{\red{A_Q}}
\newcommand{\Bq}{\red{B_Q}}
\newcommand{\Cq}{\red{C_Q}}
\newcommand{\Dq}{\red{D_Q}}

\newcommand{\Az}{\red{A_Z}}
\newcommand{\Bz}{\red{B_Z}}
\newcommand{\Cz}{\red{C_Z}}
\newcommand{\Dz}{\red{D_Z}}

\newcommand{\Acl}{{\cal A}}
\newcommand{\Bcl}{{\cal B}}
\newcommand{\Ccl}{{\cal C}}
\newcommand{\Dcl}{{\cal D}}
\newcommand{\Xcl}{{ \green{\cal X} }}

\newcommand{\B}{{\cal B}}

\newcommand{\X}{{\cal X}}

\newcommand{\Y}{\red{Y}}

\newcommand{\G}{\c{G}}

\newcommand{\kto}{\stackrel{k\to\infty}{\longrightarrow}}

\newcommand{\hl}{\\\hline}
\newcommand{\hld}{\\\hline\rule[1ex]{0ex}{1.5ex}}

\newcommand{\eql}[2]{\begin{equation}\label{#1}#2\end{equation}}
\renewcommand{\l}[1]{\label{#1}}
\let\ri\r
\renewcommand{\AA}{\ri{A}}

\renewcommand{\r}[1]{(\ref{#1})}

\newcommand{\rwc}{\rho_{\rm wc}}
\newcommand{\rro}{\rho_{\rm rs}}
\newcommand{\rec}{\rho_{\rm ec}}

\newcommand{\ch}{\check}

\newcommand{\sel}{J}

\newcommand{\ka}{\kappa}
\newcommand{\na}{\nabla}

\renewcommand{\Ac}{\red{A_a}}
\renewcommand{\Bc}{\red{B_a}}
\renewcommand{\Cc}{\red{C_a}}
\renewcommand{\Dc}{\red{D_a}}

\newcommand{\Aco}{\red{A_c}}
\newcommand{\Bco}{\red{B_c}}
\newcommand{\Cco}{\red{C_c}}
\newcommand{\Dco}{\red{D_c}}

\newcommand{\Aalp}{\Ac}
\newcommand{\Balp}{\Bc}
\newcommand{\Calp}{\Cc}
\newcommand{\Dalp}{\Dc}

\newcommand{\Aal}{{\cal A}}
\newcommand{\Bal}{{\cal B}}
\newcommand{\Cal}{{\cal C}}
\newcommand{\Dal}{{\cal D}}

\newcommand{\Gcl}{\c{G}}
\newcommand{\Xc}{\red{\c X}}

\newcommand{\At}{\t{A}}
\newcommand{\Bt}{\t{B}}
\newcommand{\Ct}{\t{C}}
\newcommand{\Dt}{\t{D}}
\newcommand{\Et}{\t{E}}
\newcommand{\Ft}{\t{F}}
\newcommand{\Mt}{\t{M}}
\newcommand{\Lt}{\t{L}}

\newcommand{\Gclt}{\t{\Gcl}}

\newcommand{\Acltf}{\t{\Acl}}
\newcommand{\Bcltf}{\t{\Bcl}}
\newcommand{\Ccltf}{\t{\Ccl}}
\newcommand{\Dcltf}{\t{\Dcl}}

\renewcommand{\Xcl}{\red{\c{X}}}

\newcommand{\cS}{\c S_{m,L}}
\newcommand{\cSs}{\c S_{0,L-m}}

\renewcommand{\le}{l}
\newcommand{\len}{l}
\newcommand{\Af}{A_{\rm f}}
\newcommand{\Bf}{B_{\rm f}}
\newcommand{\Cf}{C_{\rm f}}
\newcommand{\Df}{D_{\rm f}}
\renewcommand{\La}{\red{\Lambda}}
\newcommand{\Laopt}{\red{\Lambda^{\!*}}}
\newcommand{\Las}{\red{\bm{\Lambda}_\rho}}
\newcommand{\Lasr}{\red{\bm{\Lambda}^{r}_{\rho}}}
\newcommand{\Lai}[1]{\red{\Lambda_{#1}}}
\newcommand{\LaiT}[1]{\red{\Lambda_{#1}^T}}
\newcommand{\lai}[1]{\red{\lambda_{#1}}}

\renewcommand{\Q}{\red{Q}}
\renewcommand{\Z}{\red{Z}}

\newcommand{\Id}{I_d}
\newcommand{\z}{\text{z}}

\newcommand\ce[1]{\textcolor{blue}{#1}}

\title{Convex Synthesis of Accelerated Gradient Algorithms}

\author{Carsten Scherer\thanks{
Department of Mathematics,
University of Stuttgart, Germany, {\tt\small carsten.scherer@imng.uni-stuttgart.de}.
Funded by Deutsche Forschungsgemeinschaft (DFG, German Research Foundation) under Germany's Excellence Strategy - EXC 2075 - 390740016. We acknowledge the support by the Stuttgart Center for Simulation Science (SimTech).}
\and Christian Ebenbauer\thanks{
Institute for Systems Theory and Automatic Control, University of Stuttgart, Germany,
{\tt\small ce@ist.uni-stuttgart.de.}}%
}

\begin{document}

\maketitle
\thispagestyle{empty}
\pagestyle{empty}

\begin{abstract}
We present a convex solution for the
design of generalized accelerated gradient algorithms
for strongly convex objective functions with Lipschitz continuous gradients.
We utilize integral quadratic constraints and the Youla parameterization
from robust control theory to formulate a solution of the algorithm design problem as a convex semi-definite program.
We establish explicit formulas for the optimal convergence
rates and extend the proposed synthesis solution to extremum control problems.
\end{abstract}

\section{Introduction}

Accelerated gradient algorithms, also refereed to as
momentum methods, are considered to be among the
most widely used optimization algorithms.
These methods are applied e.g.\ in control
or artificial intelligence to train neural networks or
to solve online optimization problems arising from
receding horizon decision making.

From a control and dynamical system perspective,
accelerated algorithms can be viewed as a linear time-invariant discrete-time system in feedback with the
gradient of the to-be-minimized function as a nonlinearity
\cite{Pol87,WanEli11,DueEbe12,LesRec16,HuLes17}.
This perspective provides an immediate link to the
so-called absolute stability or Lur'e problem in control
and offers the possibility to apply advanced
tools from robust control for the analysis and design
of accelerated gradient algorithms.

It has been shown, e.g.,  
that the concept of
integral quadratic constraints and so-called Zames-Falb
multipliers allow to recover the well-known bounds for
the convergence rates of Nesterov’s celebrated
acceleration algorithm by semi-definite programming \cite{LesRec16}.
Moreover, by tuning the parameters of Nesterov’s
algorithm, these bounds can be improved to get the so-called triple momentum algorithm \cite{ScoFre18}.

A more challenging task than the analysis of
given algorithms is the design of novel algorithms
with the help of convex optimization. In light
of the relation to absolute stability
and Lur'e problems, algorithm design falls into the
area of robust feedback controller synthesis. Some
recent works have addressed the synthesis problem (e.g.\ \cite{LesSei20,MicSch20,GraEbe20,TayDro21}).
However, so far it has been an open problem to
formulate the general accelerated gradient algorithm design
problem as a genuine convex
optimization problem. In fact, the aligned question of designing robust controllers
by a convex search over the controller parameters and the multipliers to certify stability
is as well a long-standing open problem in its full generality.

In this paper, we present a convex
solution for a general accelerated gradient algorithm synthesis problem based on semi-definite programming.
Specifically, the main contributions  are as follows.
In Section~\ref{sec2}, we reveal that a particular dynamical system structure is inherent to any convergent algorithm.
This insight allows us to formulate the algorithm design problem in terms of a robust feedback controller
synthesis problem in Section \ref{Ssyn}.
We then show in Section \ref{sec_convexification}  how the special structure of the  system can be exploited to convexify the common search for the algorithm parameters and the dynamic Zames-Falb multipliers which certify convergence.
Our approach permits to derive explicit formulas for the optimal convergence rate that is achievable by synthesis,
in analogy to the analysis results for Nesterov’s algorithm in \cite{SafJos18}.
In this fashion, we are able to prove that the convergence rate of the triple momentum algorithm is indeed optimal if using the class of causal Zames-Falb multipliers to assure convergence.

Another key feature of our approach is its flexibility. We reveal in Section \ref{sec5} that it extends to extremum control \cite{AstWit89}, in which the goal is to drive the output of a dynamical system
to some steady-state condition at which a given cost function is minimized or maximized.
In particular, we establish a fully convex synthesis approach to design extremum controllers with optimal convergence properties,
even if the cost functions are structured.


Since we believe that the results in this paper are of interest to both the areas of control and optimization, we have written several sections in a tutorial fashion so that the results are accessible without a special background in robust control theory.

\section{Algorithm analysis and design}\label{Salg}
\label{sec2}

\subsection{Systems and Algorithms}

Let $\cS$ be the class of all $C^1$-functions $f:\R^d\to\R$ that
are strongly convex with parameter $m>0$ and whose gradient is Lipschitz with constant $L>m$, i.e.,
\eqn{[\na f(x)-\na f(y)]^T(x-y)\geq m\|x-y\|^2\te{and}
\|\na f(x)-\na f(y)\|\leq L\|x-y\|}
for all $x,y\in\R^d$. We denote by $\cS^0$ the set of $f\in\cS$ with $\nabla f(0)=0$.
Any $f\in \cS$ admits a unique global minimizer
$
z_*=\argmin_{z\in\R^d} f(z)\in\R^d
$
which is the solution of the equation $\na f(x)=0$.
It is well-known that the sequence defined by the gradient descent algorithm
\eql{gra}{
z_{k+1}=z_k-\al\na f(z_k)
}
for a fixed step-size $\al\in(0,2/L)$
converges to $z_*$ linearly, i.e., there exists constants $K$ and
$\rho\in( 0,1)$ such that
$\|z_k-z_*\|\leq K\rho^k\|z_0-z_*\|$ holds for all $z_0\in\R^n$ and $k\in\N_0$.
The worst-case convergence rate is defined as the infimal $\rho\in(0,1)$ for which there exists some $K$
such that linear convergence holds for all $f\in\cS$. This value depends on the algorithm parameter $\al$
and is denoted by $\rwc(\al)$. Determining upper bounds on $\rwc(\al)$
and finding an optimal choice for the algorithm parameter $\al$ which minimizes $\rwc(\al)$ has a long history in optimization theory \cite{Nes18}.

From the perspective of control, \r{gra} simply defines a nonlinear discrete-time dynamical system. Then
$k\in\N_0$ denotes a time-instant and the sequence $(z_k)_{k\in\N_0}$ is the solution (state-trajectory) of the system. Moreover, $z_*$ just constitutes a constant trajectory of \r{gra} and is, therefore,  called an equilibrium (a fixed point) thereof. Linear convergence with rate $\rho$ means that $z_*$  is globally exponentially stable with  rate $\rho$. The worst-case convergence rate is defined by considering the whole family of systems parameterized by $\na f$ for $f\in\cS$. It is common in control that such a family of systems is interpreted as a single  so-called uncertain dynamical system with an uncertainty $\na f\in \nabla\cS$.  Also in this field there is a long tradition in estimating $\rwc(\al)$, which is termed robust stability {\em analysis}. Finding a parameter $\al$ which minimizes $\rwc(\al)$ or a tight upper bound thereof is then called robust stability {\em synthesis}.

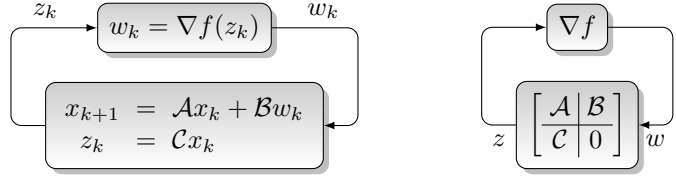
\begin{figure}[t!]
\cen{
\begin{tikzpicture}[xscale=1,yscale=1]
\node[sy2] (g) at (0,0)  {$\arr{ccl}{x_{k+1}&=&\Aal x_k+\Bal w_k\\z_k&=&\Cal x_k}$};
    \tio{Gi1}{g}{east}{1/2};
    \tio{Go1}{g}{west}{1/2};
\node[sy0] (d) at (0,1.3)  {$w_k=\na f(z_k)$};
\draw[->] (d) -| node[pos=.3,above]{$w_k$} ($(Gi1) + (1*\dl, 0)$) -- (Gi1);
\draw[->] (Go1) -- ($(Go1) - (1*\dl, 0)$) |- node[pos=.7,above]{$z_k$} (d);
\end{tikzpicture}\hspace{10ex}
\begin{tikzpicture}[xscale=1,yscale=1]
\node[sy2] (g) at (0,0)  {$\mas{c|c}{\Aal&\Bal\hl \Cal&0}$};
    \tio{Gi1}{g}{east}{1/2};
    \tio{Go1}{g}{west}{1/2};
\node[sy0] (d) at (0,1.3)  {$\na f$};
\draw[->] (d) -| ($(Gi1) + (1*\dl, 0)$) -- node[below]{$w$} (Gi1);
\draw[->] (Go1) -- node[below]{$z$} ($(Go1) - (1*\dl, 0)$) |-  (d);
\end{tikzpicture}
\caption{Feedback interconnection.}
\label{fig1}
}
\end{figure}

In robust control, a particularly useful step is to separate the description of the known parts of the algorithm from the uncertainty $\na f$. This just means to introduce the auxiliary signals $x_k:=z_k$ and
\equ{\l{del}w_k:=\na f(z_k),}
which allows us to rewrite \r{gra} as \r{del} together with
\eql{sys}{
\mat{cc}{x_{k+1}\\z_k}=\mat{cc}{\Aal&\Bal\\\Cal&0_d}
\mat{cc}{x_{k}\\w_k}
}
where $\Aal=I_d$, $\Bal=-\al I_d$, $\Cal=I_d$,  and $I_d$\,/\,$0_d$ denote the identity/zero matrix in $\R^{d\times d}$, respectively.
By itself, \r{sys} defines a linear time-invariant dynamical system that maps an initial condition $x_0\in\R^d$ and some input sequence $w=(w_k)_{k\in\N_0}$ through the recursion \r{sys} into the output sequence $z=(z_k)_{k\in\N_0}$.
The relation \r{del} alone is viewed as a static (nonlinear) system which maps the signal $z$ into $w$. Considering \r{del}-\r{sys} together means that the output (input) of \r{sys} is set equal to the input (output) of \r{del}. In control, this constitutes the feedback interconnection of \r{del} and \r{sys}
and motivates to visualize this feedback loop in an intuitive fashion by the block-diagram on the left in Fig.~\ref{fig1}. In other words, the system \r{sys} involves the algorithm parameters, while the feedback interconnection of \r{del} and \r{sys} constitutes the algorithm itself in order to compute $z_*$ for a particular instance of $f\in\cS$. Exactly the same interconnection represents an algorithm with a variable step-size $\al_k$ if just replacing $\Bal=-\al I_d$ in \r{sys} with $\Bal_k=-\al_k I_d$,
which turns \r{sys} into a linear time-varying system.

Now consider \r{sys} with general matrices $\Aal\in\R^{n\times n}$, $\Bal\in\R^{n\times d}$, $\Cal\in\R^{d\times n}$. Then the interconnection of \r{del} and \r{sys} takes the initial condition $x_0\in\R^n$ as its input and generates the unique state- and output-responses $(x_k)_{k\in\N_0}$
and $(z_k)_{k\in\N_0}$  through the recursion
\eql{usys}{x_{k+1}=\Aal x_k +\Bal \na f(\Cal x_k),\ \ z_k=\Cal x_k.}
The main goal of this work is to determine
matrices $\Aal, \Bal,\Cal$, if existing,
by a  semi-definite program
such that the algorithm \r{usys}
achieves a given convergence rate $\rho \in (0,1)$
for given $m,L$ and all objective functions $f\in\cS$.
We work with an operator interpretation of \r{sys} with general $\Aal\in\R^{n\times n}$, $\Bal\in\R^{n\times m}$, $\Cal\in\R^{k\times n}$ and replacing $0_d$ by $\Dal\in\R^{k\times m}$. Moreover, we  denote by $\le^n$ the vector space of signals $x:\N_0\to\R^n$, while $l_2^n$ is the subspace of all square summable sequences equipped with the inner product $\ab{x,y}_2:=\sum_{k=0}^\infty x_k^Ty_k$ and the norm $\|x\|_2:=\sqrt{\ab{x,x}_2}$
for $x,y\in l_2^n$. For $w\in \le^m$ and $x_0\in\R^n$, the recursion
\eql{rec}{
\mat{cc}{x_{k+1}\\z_k}=\mat{cc}{\Aal&\Bal\\\Cal&\Dal}
\mat{cc}{x_{k}\\w_k}
}
defines unique state- and output-responses $x\in \le^n$ and $z\in\le^d$, respectively; for a fixed $x_0$ (often taken to be $0$), the resulting affine (linear) operator is denoted as
\eql{sysm}{
z=\mas{c|c}{\Aal&\Bal\hl\Cal&\Dal}w
}
where we suppress the dependence on $x_0$. All throughout this paper,
we reserve square brackets to denote the input-output operator \r{sysm} defined by the recursion \r{rec};
the partition lines in \r{sysm} are always displayed to separate $\Aal$ from the other matrix blocks in \r{rec}.

Consequently \r{del}-\r{sys}, \r{usys}, and
\eql{nl}{w=\nabla f(z),\ \ z=\mas{c|c}{\Aal&\Bal\hl\Cal&0}w}
express one and the same interconnection as
depicted by the block-diagram in Fig.~\ref{fig1}.
Analyzing the convergence properties of a general algorithm \r{usys} then boils down to analyzing
the stability properties of the feedback system \r{nl}.  From now on we represent algorithms interchangeably by \r{sys} or \r{sysm}. Algorithm convergence means that the signal $z$ in \r{nl} converges to the minimizer of $f$ for any initial condition $x_0\in\R^n$.

The analysis of stability of feedback interconnections constitutes one of the fundamental questions studied in control since its beginnings, with many traditional ideas nicely collected in the classical textbook \cite{DesVid75}.
Polyak was among the first to clearly emphasize the above sketched tight link between the two areas \cite{Pol87}, see also \cite{WanEli11,DueEbe12,LesRec16,HuLes17}. By arguing with an analogy to mechanical systems, he suggested to replace \r{gra} by the heavy-ball method which includes a damping or momentum term as in $z_{k+1}=z_k-\al \na f(z_k)+\be(z_k-z_{k-1})$ for some parameters $\be\in[0,1)$ and $\al\in(0,2(1+\beta)/L)$. It is easy to check that the corresponding algorithm is \r{nl} for the system matrices
\eql{tri}{
\mat{c|c}{\Aal&\Bal\hl\Cal&0}=\mat{cc|c}{(1+\be)I_d&-\be I_d&-\al I_d\\I_d&0&0\hl (1+\ga)I_d&-\ga I_d&0}
}
with $\ga=0$. 
In \cite{Pol87} it is shown that the convergence rate is considerably improved over gradient descent, at the cost of sacrificing global algorithm convergence \cite{LesRec16}. Nesterov's celebrated accelerated gradient decent algorithm corresponds the choice $\ga=\be$ in \r{tri} with guaranteed global and fast convergence \cite{Nes18}. The more recently proposed triple momentum algorithm \cite{ScoFre18} relies on different values of the three parameters in \r{tri} with the best-known convergence rate to date.

Let us conclude this section by recalling some basic notions for general linear systems \r{rec} or \r{sysm}.
With an invertible matrix $T\in\R^{n\times n}$, a state-coordinate change for \r{rec} is defined by $\xi_k:=Tx_k$. It is easily seen that this transforms the quadrupel $(\Aal,\Bal,\Cal,\Dal)$ into $(T^{-1}\Aal T,T^{-1}\Bal,\Cal T,\Dal)$.
For a fixed input signal $w$ and the initial conditions $x_0$ and $\xi_0=Tx_0$, respectively,
one can check that the output trajectories of the original and the transformed systems are identical;
this is compactly expressed through
\eql{sys1}{
\mas{c|c}{\Aal&\Bal\hl\Cal&\Dal}=
\mas{c|c}{T^{-1}\Aal T&T^{-1}\Bal\hl\Cal T&\Dal}.
}
Furthermore, the series interconnection of two systems
\eql{sys12}{
y_1=\mas{c|c}{\Aal_1&\Bal_1\hl \Cal_1&\Dal_1}u_1\te{and}
y_2=\mas{c|c}{\Aal_2&\Bal_2\hl \Cal_2&\Dal_2}u_2
}
is defined by using the output signal of the second as an input to the first, which is reflected by $u_1=y_2$ (and requires that the  signal dimensions match). This is nothing but the composition of the two respective maps,
which is as usual denoted as an operator product. It is elementary to verify that the series interconnection of
the two systems \r{sys12}
can be described by
\equ{\label{prod}\arc
\mas{cc|c}{\Aal_1&\Bal_1\Cal_2&\Bal_1\Dal_2\\0&\Aal_2&\Bal_2\hl \Cal_1&\Dal_1\Cal_2&\Dal_1\Dal_2}\text{ or }
\mas{cc|c}{
\Aal_2    &0   &\Bal_2    \\
\Bal_1\Cal_2 &\Aal_1 &\Bal_1\Dal_2 \hl
\Dal_1\Cal_2 &\Cal_1 &\Dal_1\Dal_2 }.
}

In case of identical dimensions of the input and output signals in \r{sys12}, the sum of the two maps is the so-called parallel interconnection given by
\eql{sum}{\arraycolsep.4ex
\mas{c|c}{\Aal_1&\Bal_1\hl \Cal_1&\Dal_1}\!+\!\mas{c|c}{\Aal_2&\Bal_2\hl \Cal_2&\Dal_2}=
\mas{cc|c}{\Aal_1&0&\Bal_1\\0&\Aal_2&\Bal_2\hl \Cal_1&\Cal_2&\Dal_1\!+\!\Dal_2}.}
Further, if $\Dal$ is invertible, the map \r{sysm} is invertible and its inverse can be represented with
$$
\mas{c|c}{\Aal-\Bal\Dal^{-1}\Cal&\Bal\Dal^{-1}\hl -\Dal^{-1}\Cal&\Dal^{-1}}.
$$
The system \r{rec} is called stable if $\Aal$ is a Schur matrix, i.e., all its eigenvalues are in absolute value strictly smaller than one. Moreover, \r{rec} or the pair $(\Aal,\Bal)$ is stabilizable if there exists a matrix $M$ such that $\Aal+\Bal M$ is Schur; similarly, \r{rec} or $(\Aal,\Cal)$ is detectable if there exists $L$ such that $\Aal+L\Cal$ is Schur.

\subsection{Algorithm Structure and Convergence}\label{Sstr}

Let us now get back to the algorithm \r{nl}.
First, we settle that convergence enforces an important
structural constraint on the parameters $\Aal$, $\Bal$, $\Cal$.
We start by stressing that there is no benefit in choosing systems \r{sys} which are not detectable.
Indeed, suppose $(\Aal,\Cal)$ in \r{nl} is not detectable. We then follow
\cite[Sec. 3.3]{ZhoDoy96} and perform a state-coordinate change to obtain
$$
\mat{c|c}{\Aal&\Bal\hl \Cal&0}=\mat{cc|c}{\Aal_1&\Aal_{12}&\Bal_1\\0&\Aal_2&\Bal_2\hl 0&\Cal_2&0}
$$
where $(\Aal_2,\Cal_2)$ is detectable. Due to the block structure of the matrices,
the set of $z$-trajectories of \r{usys} and of
$\xi_{k+1}=\Aal_2 \xi_k +\Bal_2 \na f(\Cal_2 \xi_k),\ \ z_k=\Cal_2 \xi_k$
are obviously identical. W.l.o.g. we can hence replace the non-detectable system in \r{nl}
by the one described with the triple $(\Aal_2,\Bal_2,\Cal_2)$ which is detectable.

Note that the least convergence requirement for the algorithm \r{del}-\r{sys} is
\eql{conv}{\lim_{k\to\infty}z_k=z_*\te{and}\lim_{k\to\infty}w_k=0}
(for any initial condition $x_0\in\R^n$ and some $z_*\in\R^d$).
If $(\Aal,\Cal)$ is detectable and we take $L$ such that $\Aal+L\Cal$ is Schur, we infer
$x_{k+1}=(\Aal+L\Cal)x_k+\Bal w_k-Lz_k$. Therefore, \r{conv} also implies the convergence of the state-trajectory $x_k$ to some $x_*$ for $k\to\infty$ with
\eql{equ}{x_*=\Aal x_*\te{and}z_*=\Cal x_*.}
Most importantly, we now show that \r{conv} enforces the following special structure of the algorithm parameters.
\theo{\label{Tnec}Let $(\Aal,\Cal)$ be detectable. If all trajectories of \r{usys} satisfy \r{conv} for all quadratic functions $f\in\cS$ and all $x_0\in\R^n$, then $\Aal+\Bal m\Cal$ is Schur and there exist 
$\Aalp,\Balp,\Calp,\Dalp$ such that
\eql{imp}{
\mas{c|c}{\Aal&\Bal\hl \Cal&0}=
\mas{c|c}{\Aalp&\Balp\hl\\[-2.5ex] \Calp&\Dalp}
\mas{c|c}{I_d&I_d\hl I_d&0}.
}
If $(\Aal,\Cal)$ has the structure induced by
\r{imp}, then \r{equ} has a unique solution $x_*\in\R^n$ for every $z_*\in\R^d$.
}

Before entering the proof, let us interpret the structural property \r{imp} in the state-space. By \r{prod}, it implies that there exists a state-coordinate change of \r{sys} after which the algorithm \r{del}-\r{sys} reads
\eql{alg}{
\arc
\mat{c}{\xi_{k+1}\\\eta_{k+1}\hl z_k}=
\mat{cc|c}{
I_d   &0     &I_d \\
\Balp &\Aalp &0   \hl
\Dalp &\Calp &0   }
\mat{c}{\xi_{k}\\\eta_{k}\hl w_k},\ w_k=\nabla f(z_k).
}
\proof{We start by observing that, given $z_*$, there is at most one vector $x_*$ satisfying \r{equ} since
\eql{frk}{\rk\mat{c}{\Aal-I\\\Cal}=n.}
Indeed, because $(\Aal,\Cal)$ is detectable, we can take $L$ such that $\Aal+L\Cal$ is Schur; then
$(\Aal-I)x=0$, $\Cal x=0$ imply $(\Aal+L\Cal-I)x=0$ and thus $x=0$, because $1$ is no eigenvalue of $\Aal+L\Cal$.

Now take $f\in\cS$ as $f(z)=\frac{1}{2}z^T(m\Id)z-b^Tz$
with 
any $b\in\R^d$. Then the trajectories of
\r{del}-\r{sys} satisfy
\eql{lusys}{
x_{k+1}=(\Aal +\Bal m\Cal )x_k-\Bal b,\ \ w_k=m\Cal x_k-b.
}
Due to \r{conv} we have argued above that $x_k\kto x_*$ with $x_*$ satisfying \r{equ}; we also infer $m\Cal x_*=b$.

Let us first take $b=0$. By $m>0$ we get $z_*=\Cal x_*=0$.
Since $x_*=0$ satisfies \r{equ}, we infer (by uniqueness) that
$x_k\kto 0$ for all trajectories of \r{lusys}. This implies that $\Aal+\Bal m\Cal$ is Schur which proves the first statement.

Now let $b\in\R^d$ be general. We then get $b=m\Cal x_*=m\Cal(\Aal +\Bal m\Cal -I_n)^{-1}\B b$
for all $b\in\R^d$ and, therefore, $I_d-m\Cal(\Aal+\Bal m\Cal-I_n)^{-1}\Bal =0$. With a Schur complement argument \cite{HorJoh85}, this implies
\eqn{
\rk\mat{cc}{\Aal +\Bal m\Cal -I_n&\Bal \\m\Cal &I_d}=n.
}
Yet another Schur complement argument shows
\eql{h1}{
\rk(\Aal -I_n)=n-d.
}
Now choose $T_1\in\R^{n\times d}$ with full column rank and $T_1^T(\Aal-I_n)=0$. Then $T_1^T\Bal \in\R^{d\times d}$ is invertible;  otherwise there exists $x\neq 0$ with $x^TT_1^T\Bal =0$ and thus
$(T_1x)^T(\Aal +\Bal m\Cal)=x^TT_1^T\Aal =x^TT_1^T=(T_1x)^T$; because $T_1x\neq 0$, we infer
that $1$ is an eigenvalue of $\Aal +\Bal m\Cal$, which is a contradiction since the latter matrix is Schur. We can hence choose $T_1$ to also satisfy
$T_1^T\Bal=I_d$ and pick $T_2\in\R^{n\times (n-d)}$ with
$T_2^T\Bal=0$ such that $T=(T_1\ T_2)$ is invertible. We get
$
T^T\Aal =\mat{cc}{I_d&0\\\Balp&\Aalp}T^T\te{and}T^T\Bal=\mat{cc}{I_d\\0}
$
for suitable matrices $\Aalp\in\R^{(n-d)\times(n-d)}$, $\Balp\in\R^{(n-d)\times d}$.
With $\mat{cc}{\Dalp&\Calp}:=\Cal T^{-T}$ we infer \r{imp} since
\eql{h2}{
\mat{c|c}{T^T\Aal T^{-T}&T^T\Bal \hld \Cal T^{-T}&0}=
\mat{cc|c}{
I_d   &0     &I_d \\
\Balp &\Aalp &0   \hl
\Dalp &\Calp &0   }.
}

For the system matrices in \r{alg}, we finally note that \r{equ} is equivalent to
\eql{h3}{
\mat{ccc}{\Balp&\Aalp-I\\\Dalp&\Calp}x_*=\mat{c}{0\\z_*}
}
and that the matrix ín \r{h3} is square. Since \r{h3} has at most one
solution $x_*$ (as shown at the beginning of the proof), we infer that the matrix in \r{h3}
is actually invertible, which proves the last statement.\epro}



We are now in the position to introduce the precise definition of algorithm convergence with rate $\rho$.
\defi{
Let the system in \r{nl} be detectable and admit the structure \r{imp}.
For $\rho\in(0,1)$, algorithm \r{nl} achieves $\rho$-convergence
(for the class $\cS$) if
there exists some $K\geq 0$ such that
\eql{conrho}{\|z_k-z_*\|\leq K\rho^k\|x_0-x_*\|\te{for all}k\in\N_0,}
for any $f\in\cS$ with minimizer $z_*\in\R^d$, any $x_0\in\R^n$
and any $x_*\in\R^n$ satisfying \r{equ}.

The infimum of all $\rho\in(0,1)$ such that \r{nl} achieves $\rho$-convergence is the
algorithm convergence rate and denoted as $\rwc$ (with $\rwc:=\infty$ if no such $\rho$ exists).}

Note that $\rho$-convergence is invariant under a state-coordinate change of \r{sys}.
Moreover, $\rho$-convergence implies but is stronger than the convergence property \r{conv} for all trajectories of \r{nl} with any $f\in\cS$.

\theo{\label{Tsuf}
Let the system in \r{nl} be detectable and admit the structure \r{imp}.
If \r{nl} achieves $\rho$-convergence
for the class $\cS^0$, then it achieves $\rho$-convergence for the full class $\c{S}_{m,L}$ as well.
}

\proof{By assumption, there exists some $K\geq 0$ such that all trajectories of
\eql{sys0}{
{\arc
\mat{c}{\bar x_{k+1}\\ \bar z_k}}=
\mat{cc}{\Aal&\Bal\\\Cal&0}
\mat{c}{\bar x_{k}\\ \bar w_k},\ \bar w_k=\nabla \bar f(\bar z_k)
}
for any $\bar f\in\cS^0$ satisfy
\eql{int1}{\|\bar z_k\|\leq K\rho^k\|\bar x_0\|\te{for all}k\in\N_0.}
Now take $f\in\cS$ with minimizer $z_*$, any $x_0\in\R^n$ and consider \r{alg}.
By Theorem \ref{Tnec}, \r{equ} has a unique solution $x_*\in\R^n$ which clearly satisfies
$
\mat{c}{x_*\\z_*}=\mat{cc}{\Aal&\Bal\\\Cal&0}\mat{c}{x_*\\ 0}.
$
Define $\bar x_k:=x_k-x_*$, $\bar z_k:=z_k-z_*$, $\bar w_k:=w_k$ and $\bar f:=f(\bullet+z_*)$. By linearity,
this yields a trajectory of \r{sys0}. Since $\bar f\in\cS^0$, we infer that \r{int1} is valid. This is clearly nothing but \r{conrho} as was to be shown.
\mbox{}\hfill\epro
}

In summary, Theorem \ref{Tnec} reveals that algorithm convergence requires that the related linear system ``contains'' a model of the so-called discrete time integrator
\eql{int}{
\mat{c}{\eta_{k+1}\\y_k}=\mat{cc}{I_d&I_d\\I_d&0}\mat{c}{\eta_k\\u_k}
}
as a right factor. Conversely, by Theorem \ref{Tsuf},
if the algorithm parameters "contain" the integrator \r{int}, $\rho$-convergence can be induced from $\rho$-convergence for $f\in\cS^0$ with a minimizer located at the origin. From a control theory perspective, this is reminiscent of the so-called internal model principle \cite{Won85a}.

\subsection{Robust Stability Analysis and $\rho$-Convergence}\label{Sana}

From now on we assume that the system in \r{nl} is detectable and admits the structure \r{imp}.
The next goal is to relate the question of bounding the algorithm convergence rate $\rwc$ to a robust stability analysis problem.
In view of Theorems \ref{Tnec} and \ref{Tsuf}, it suffices to confine the discussion to the class $\cS^0$. 
We also  map $\cS^0$ bijectively onto $\cSs^0$ through $f\mapsto g$ where
$g(z):=f(z)-\frac{1}{2}z^T(m\Id)z\te{for}z\in\R^d.$
Then, \r{nl} clearly just is the interconnection of
\r{sys} with $w=\nabla g(z)+mz$. With the transformation
\eql{sit0}{
\mat{c}{\bar z\\\bar w}:=\mat{cc}{\Id&0\\-m\Id&\Id}\mat{c}{z\\w},
}
this interconnection can be as well expressed by
\eql{cll1}{\arc
\mat{c}{\bar x_{k+1}\\\bar z_k}=\mat{cc}{\Aal+\Bal m\Cal&\Bal\\ \Cal&0}\mat{c}{\bar x_k\\\bar w_k},\ \bar w_k=\nabla g(\bar z_k).
}
Then $\rwc$ is just equal to the convergence rate of \r{cll1} for
the class $\cSs^0$.
Next, for $\rho\in(0,1)$, we follow \cite{DesVid75} and introduce the signal weighting mapping
$
\rho_+:\le^n\to\le^n,\ \ x\mapsto \rho_+(x)=(\rho^kx_k)_{k\in\N_0}
$
which is bijective. Then
\eql{sit1}{
\ch x:=\rho_+^{-1}(\bar x),\ \ \ch w=\rho_+^{-1}(\bar w)\te{and}
\ch z=\rho_+^{-1}(\bar z)}
transform \r{cll1} into
\eql{cll}{
\mat{c}{\ch x_{k+1}\\\ch z_k}=\mat{cc}{\rho^{-1}(\Aal+\Bal m\Cal)&\rho^{-1}\Bal\\ \Cal&0}\mat{c}{\ch x_k\\\ch w_k},\ \
\ch w_k=\rho^{-k}\nabla g(\rho^k \ch z_k).
}
These steps permit us to relate $\rwc$ to a standard robust stability margin for the map $\ch x_0\mapsto \ch z$ defined by \r{cll}.

\lemm{\label{Lrs}Let $\rro$ be the infimal $\rho\in(0,1)$ for which there is a $K\geq 0$ such that all trajectories of \r{cll} with $g\in \cSs^0$ satisfy $\|\ch z\|_2\leq K\|\ch x_0\|$.
Then the convergence rate $\rwc$ of algorithm \r{nl} is equal to $\rro$.
}

\proof{To show $\rwc\leq\rro$ we can assume $\rro<\infty$.
Let $\rho\in(\rro,1)$ and take any trajectory of \r{cll1}.
Then the $\rho_+^{-1}$-transformed signals define a trajectory of \r{cll} and we hence infer, by using the definition of $\rro$, that
$
\rho^{-k}\|\bar z_k\|=\|\ch z_k\|\leq \|\ch z\|_2\leq K\|\ch x_0\|=K\|\bar x_0\|
$
for all $k\in\N_0$. This proves $\rho$-convergence of \r{cll1} for the class $\cSs^0$ and, therefore,
$\rwc\leq\rho$. Since $\rho\in(\rro,1)$ was arbitrary, we conclude $\rwc\leq\rro$.

To see $\rro\leq\rwc$ let $\rwc<\infty$ and take $\rho\in(\rwc,1)$. Choose
some $\rho_0\in(0,1)$ with $\rho_0\rho\in(\rwc,1)$ to infer $(\rho_0\rho)$-convergence of \r{cll1}
for the class $\cSs^0$. Hence, there exists some $\bar K\geq 0$ such that all trajectories of \r{cll1} satisfy
$\|(\rho_0\rho)^{-k}\bar z_k\|\leq \bar K\|\bar x_0\|$ and thus
$\|\rho^{-k}\bar z_k\|\leq \bar K\rho_0^{k}\|\bar x_0\|$
for all $k\in\N_0$. Then any trajectory of \r{cll} can be transformed with $\rho_+$ back into one of \r{cll1} to get, with $\bar x_0=\ch x_0$, that
$
\sum_{k=0}^\infty \|\ch z_k\|^2=\sum_{k=0}^\infty \|\rho^{-k}\bar z_k\|^2\leq
\left(\sum_{k=0}^\infty\rho_0^{2k}\right)\bar K^2\|\ch x_0\|^2.
$
We conclude $\rro\leq\rho$ and, hence, $\rro\leq\rwc$. \epro
}

In summary, computing tight upper bounds on the convergence rate $\rwc$ of algorithm \r{nl} is equivalent to
determining tight upper bounds on the so-called robust stability margin $\rro$ for the interconnection \r{cll} as defined in Lemma \ref{Lrs}.

\remark{\label{rem1} If $\rro<\infty$ then
$\rho^{-1}(\Aal+\Bal m\Cal)$ is Schur for all $\rho\in(\rro,1)$. This is shown as in the first step of the proof of
Theorem \ref{Tnec}.
}


\subsection{Robust Stability Analysis and Integral Quadratic Constraints}\label{Siqc}

In this section we sketch how to compute effective bounds on the margin $\rwc=\rro$
by setting up a semi-definite program. This involves a family of so-called Zames-Falb multipliers.
These are systems
\eql{fil}{\arc
\Pi(\La):=
\mas{c|c}{\Af&\Bf\hl \Cf(\La)&\Df(\La)}:=
\mas{cccc|c}{
0      &I_d    &\cdots &0      &0        \\
\vdots &\ddots &\ddots &\vdots &\vdots   \\
0      &\cdots &0      &I_d    &0      \\
0      &\cdots &0      &0      &I_d      \hl
\Lai{\len}&\Lai{\len-1}&\cdots&\Lai{1}&\Lai{0}}
}
which are parameterized by a matrix tuple $\La$ in the set
\mun{
\Las:=\left\{\mat{cccc}{\Lai{\le}&\cdots&\Lai{1}&\Lai{0}}
\in(\R^{d\times d})^{l+1}
\st
\right.
\Lai{0}-\diag(\Lai{0})\leq 0,\ \ \Lai{i}\leq 0\te{for}i=1,\ldots,l,\\
\left. \left(\sum_{i=0}^\len \Lai{i}\rho^i\right)e\geq 0,\
e^T\left(\sum_{i=0}^\len \Lai{i}\rho^{-i}\right)\geq 0
\right\};
}
the inequalities are read elementwise and $e\in\R^d$ is the all-ones vector, while
$\diag(A)\in\R^{d\times d}$ is the diagonal matrix whose diagonal is identical to that of $A\in\R^{d\times d}$.

Note that \r{fil} is a so-called finite-impulse-response filter of length $\le$ and of dimension $d\times d$; the latter two parameters are not displayed in $\Las$ to lighten the notation. The parameters for which
the blocks in $\La\in\Las$ are diagonally repeated are collected in
\eql{Lasr}{
\Lasr\!:=\Las\cap\left\{\arc\mat{cccc}{\lai{\len}I_d&\cdots&\lai{1}I_d&\lai{0}I_d}\,|\,
\lai{i}\in\R\right\}.
}

The introduction of this family is motivated by the robust stability result in
\cite{MicSch20} for \r{cll}.
This involves the following positivity property for the nonlinearity
\eql{grag}{
\ch w=\rho_+^{-1}\nabla g(\rho_+\ch z).}

\lemm{\label{Liqc}Let $g\in\cSs^0$ and $\La\in\Lasr$. Then
\eql{pas}{\ab{\Pi(\La)\t z,\t w}_2\geq 0}
holds for all $\ch z\in l_2^d$ and the output $\ch w$ of \r{grag} with
\eql{sit}{
\mat{c}{\t z\\\t w}:=\mat{cc}{L\Id-m\Id&-\Id\\0&\Id}\mat{c}{\ch z\\\ch w}.
}
}

We emphasize that  the response $\Pi(\La)\t z$ is defined based on the state-space representation \r{fil} with
the state's initial condition taken as zero. From now on, we follow this convention
in robust control for systems expressed in operator notation as in \r{sysm}.

In systems theory, \r{pas} is a so-called passivity property for filtered versions of the input and output signals of \r{grag}; it is
also referred to as an integral quadratic constraint (IQC) \cite{MegRan97}.
The latter terminology emerges since such results are often formulated for continuous time systems, for which the $l_2$-inner product (expressed in terms of sums) is replaced by the inner product on $L_2[0,\infty)$ (involving integrals).

Guaranteeing robust stability of \r{cll} involves a related negativity condition for the linear system
\eql{sysg}{
\ch z=
\mas{c|c}{\rho^{-1}(\Aal+\Bal m\Cal)&\rho^{-1}\Bal\hl \Cal&0}
\ch w.
}
By Remark \ref{rem1}, this system needs to be stable. As a consequence, if $\ch w\in l_2^d$ is any input with finite $l_2$-norm,  the response of \r{sysg} satisfies $\ch z\in l_2^d$. The following
filtered strict negativity property then guarantees
robust stability for \r{cll} as defined in Lemma \ref{Lrs} and, thus, assures $\rro=\rwc\leq\rho$ \cite[Lemma 3, Theorem 4]{MicSch20}.
\theo{\label{Trs}
Let $\rho\in(0,1)$. Then $\rwc\leq\rho$ is assured if
$\rho^{-1}(\Aal+\Bal m\Cal)$ is Schur and if there exist $\La\in\Lasr$, $\eps>0$ such
that for any $\ch w\in l_2^d$ and the response $\ch z\in l_2^d$ of \r{sysg},
the signals \r{sit} satisfy
\eql{pass}{\ab{\Pi(\La)\t z,\t w}_2\leq -\eps\|\t w\|_2^2.}}
In view of \r{sit} and \r{sysg} and for any $\t w\in l_2^d$, the trajectories $(\t z,\t w)$ in \r{pass} can as well be associated
to the system $\t z=\tilde{\c G}\t w$ with
\eql{keysys}{
\Gclt:=\mas{c|c}{\rho^{-1}(\Aal+\Bal m\Cal)&\rho^{-1}\Bal\hl (L-m)\Cal&-I_d}.}
If $G$ denotes the series interconnection $\Pi(\La)\t{\c G}$, then
\r{pass} just reads $\ab{G\t w,\t w}_2\leq-\eps\|\t w\|_2^2$ for all $\t w \in l_2^d$,
and $G$ is also said to be strictly negative real. Theorem~\ref{Trs} just expresses that
$\rwc\leq\rho$ is guaranteed by checking that $G$ is stable and strictly negative real.

To verify these properties, we can use the following variant of the celebrated positive real lemma \cite{AndVon73}.

\lemm{\label{Lpr}Consider a system
$
z=Gw\text{ with }G=\mas{c|c}{A&B\hl C&D}
$
and $d$ inputs and outputs. Then the following statements are equivalent:
\enu{
\item $A$ is Schur and $G$ is strictly negative real (SNR): There exists some $\eps>0$ with $\ab{Gw,w}_2\leq -\eps\|w\|_2^2$ for all $w\in l_2^d$.
\item
There exists some $\red{X}\cg 0$ such that
\eql{negrea}{\arc
\mat{cc}{A&B\\I&0\hl C&D\\0&I}^T\!\!\!
\mat{cc|cc}{\red{X}&0&0&0\\0&-\red{X}&0&0\hl 0&0&0&I\\0&0&I&0}
\mat{cc}{A&B\\I&0\hl C&D\\0&I}\cl 0.
}}
}
In here, $A\cg 0$ ($A\cl 0$) means that the real matrix $A$ is symmetric and positive (negative) definite.
Lemma \ref{Lpr} allows to translate stability and strict negative realness of an operator defined by a linear system into a convex finite-dimensional feasibility constraint, which takes of the form of a linear matrix inequality (LMI) in the matrix variable $\red{X}$.

Recall that, in Theorem \ref{Trs}, this involves the series interconnection of $\Pi(\La)$ and $\tilde{\c G}$ in \r{keysys} with the
state-space description
\equ{\label{wcll}
\Pi(\La)\tilde{\c G}=
\mas{c|c}{\Acltf&\Bcltf\hld \Ccltf(\La)&\Dcltf(\La)}:=
\mas{cc|cc}{
\Af     &\Bf(L-m) \Cal               &-\Bf          \\
0       &\rho^{-1}(\Aal+\Bal m\Cal)   &\rho^{-1}\Bal \hl
\Cf(\La)&\Df(\La)(L-m)\Cal           &-\Df(\La)}.}
We observe that $\rho^{-1}(\Aal+\Bal m\Cal)$ is Schur iff this holds for $\Acltf$.
Hence, Theorem \ref{Trs} in combination with Lemma~\ref{Lpr} leads to the following result.

\coro{\label{Can}
For $\rho\in(0,1)$, the convergence rate of algorithm \r{nl} is bounded as
$\rwc\leq\rho$ if there exist $\La\in\Lasr$ and $\Xc$ that satisfy the constraints
\gat{
\label{lmia}\arc
\Xc\cg 0,\ \
\bul^T
\mat{cc|cc}{\Xc&0&0&0\\0&-\Xc&0&0\hl 0&0&0&I\\0&0&I&0}
\mat{cc}{\Acltf&\Bcltf\\I&0\hld \Ccltf(\La)&\Dcltf(\La)\\0&I}\cl 0.
}
}
In \r{lmia} and later we use ``$\bul$'' as a placeholder for the matrix on the right to save space.
For fixed $l\in\N_0$ and $\rho\in(0,1)$, we observe that the
left-hand side of \r{lmia} is affine in the variables $\Xc$ and $\La$.
Note that $\Lasr$ is as well described by LMI constraints on $\La$.
We have thus reduced the verification of $\rho$-convergence for a {\em given} algorithm
to a convex feasibility test in terms of LMIs.

However, recall that the main goal of this paper is algorithm {\em design}. For a given $\rho\in(0,1)$, this means to verify
with an LMI feasibility test whether there exists an algorithm of the form \r{usys} which achieves $\rho$-convergence;
if existing, one wishes to construct an algorithm based on some solution of the LMI.

Recall that we parameterize algorithms by the matrices $\Aalp,\Balp,\Calp,\Dalp$ in \r{alg}.
In Corollary~\ref{Can}, these matrices enter the constraints \r{lmia} via \r{wcll} in a non-linear fashion.
If using Corollary \ref{Can} for design, we end up with non-convex constraints
if viewing the algorithm matrices $\Aalp,\Balp,\Calp,\Dalp$, the multiplier parameter $\La$ and the so-called
Lyapunov matrix $\Xc$ as decision variables. As a preparation for overcoming this trouble,
we recapitulate some essential insights into controller synthesis by LMIs next. 


%

\section{Controller Synthesis with LMIs}\label{Ssyn}

Feedback control for systems described by difference equations as in this paper can be abstractly formulated in terms of a given to-be-controlled system, the so-called generalized plant, as described by
\eql{gpla}{
\mat{cc}{z\\y}=
\mas{c|ccc}{A&B_1&B\hl C_1&D_1&E\\C&F&0}
\mat{cc}{w\\u}.}
This plant has two (vector-valued) input and output signals. Here $u$ is the so-called control input with which the system is actuated, steered or manipulated. On the other hand, $y$ is the so-called measurement output, which is viewed as the available information about the system for the purpose of controlling it. A controller then takes $y$ as its input and generates the control action $u$ as its output through
\eql{gcon}{u=\mas{c|c}{\Aco&\Bco\hl \Cco&\Dco}y.}
The plant and controller form the so-called feedback interconnection, which can be expressed
(after  a simple calculation \cite[Section 2]{Sch00siam}) as
\eql{gcll}{
z=
\mas{cc|c}{
A+B\Dco  C   &B\Cco &B_1+B\Dco F  \\
\Bco C       &\Aco  &\Bco F       \hl
C_1+E\Dco C  &E\Cco &D_1+E\Dco F  }w=:\mas{c|c}{\Acl&\Bcl\hl \Ccl&\Dcl}w.}
The closed-loop system is affected by the disturbance input $w$ and responds with the controlled output $z$; these are the signals on which one imposes desired specifications which the controller
should achieve.

Foremost, controllers are required to (internally) stabilize the plant, i.e., they
need to render $\Acl$ Schur.
Next to stabilization, many desired so-called performance properties on the map $w\mapsto z$ are expressed as
$$
\ab{\mat{c}{z\\w},P\mat{c}{z\\w}}_2\leq-\eps\|w\|_2^2
$$
for all $w\in l_2^d$ (where $d$ is the number of components of $w$) and some $\eps>0$. Here $P$ is an indefinite symmetric weighting matrix that is partitioned according to the signals $z$ and $w$ with the properties
$$
P=\mat{cc}{Q&S\\S^T&R},\ \ Q\cge 0\te{and}\det(P)\neq 0.
$$
The celebrated KYP lemma (see e.g. \cite{ZhoDoy96,Ran96})
can be used to show that the controller achieves both tasks iff
there exists a 
Lyapunov matrix $\Xcl$ that satisfies
\eql{glmia}{\arc\Xcl\cg 0,\ \bul^T
\mat{cc|cc}{\Xcl&0&0&0\\0&-\Xcl&0&0\hl 0&0&Q&S\\0&0&S^T&R}
\mat{cc}{\Acl&\Bcl\\I&0\hl \Ccl&\Dcl\\0&I}\cl 0.}
Analyzing the desired properties of a fixed controller thus boils down to this convex feasibility test in $\Xcl$.


If synthesizing a controller, we view the parameters of \r{gcon} as additional decision variables. However, \r{glmia} does not impose a convex constraint on both $\Xcl$ and the controller matrices. Despite this trouble, the existence of a controller that achieves    \r{glmia} can still be equivalently expressed as convex constraints.

\theo{\label{Tsyn}Let $U$ and $V$ be matrices whose columns form a basis of $\arc\ker\mat{cc}{C&F}$ and $\arc\ker\mat{cc}{B^T&E^T}$, respectively. Then there exist a controller \r{gcon}
and an $\Xcl$ such that the closed-loop system \r{gcll} satisfies
\r{glmia} iff there exist symmetric matrices $\red{X}$, $\red{Y}$ with
\eql{lmis1}{\arc
\mat{cc}{\red{X}&I\\I&\red{Y}}\cg 0,\
\bul^T
\mat{cc|cc}{
\red{X}&0        &0 &0 \\
0      &-\red{X} &0 &0 \hl
0      &0        &Q &S \\
0      &0        &S^T &R }
\mat{cc}{A&B_1\\I&0\hl C_1&D_1\\0&I}U\cl 0,\
\bul^T
\mat{cc|cc}{\red{Y}&0&0&0\\0&-\red{Y}&0&0\hld 0&0&\t Q&\t S\\0&0&\t S^T&\t R}
\mat{cc}{I&0\\-A^T&-C_1^T\hl 0&I\\ -B_1^T&-D_1^T}V\cg 0;
}
here $\t Q$, $\t R$, $\t S$ denote the blocks of the inverse $P^{-1}$.
}

Once the LMIs \r{lmis1} are feasible, a constructive procedure to compute the controller matrices is found in \cite{Gah94}, and the dimension of the resulting state-matrix $\Aco$ equals that of $A$. 

This result essentially appeared in the  seminal work \cite{GahApk94,IwaSke94}
for $Q=I$, $R=-\ga^2I$, $S=0$ related to the so-called $H_\infty$-control problem.
The extensions to general performance indices have been suggested in \cite{MasOha98,SchGah97},
while the current paper is aligned in notation with the exposition in
\cite{Sch00siam,SchWei99}. Here we only use
\eql{QRS}{Q=R=0\te{and}S=I.}
\remark{\label{rem2}
To ensure the existence of a stabilizing controller for \r{gpla},
one should verify that $(A,B)$ is stabilizable and $(A,C)$ is detectable \cite[Section 17.1]{ZhoDoy96}.
}

\section{Convexification of Algorithm Synthesis}
\label{sec_convexification}

\subsection{Algorithm Design by Controller Synthesis}\label{Ssyna}

It is now a natural idea to exploit the general controller synthesis framework in Section \ref{Ssyn} for
algorithm design based on Corollary \ref{Can}.

In order to match \r{lmia} with \r{glmia}, we choose \r{QRS}
and express the state-space description of \r{wcll} as the interconnection of a suitable
plant \r{gpla} in feedback with a controller which is determined through the algorithm parameters. 
In fact, a trivial computation shows that the system \r{keysys} for
$(\Aal,\Bal,\Cal)$ from \r{alg} is obtained as in \r{gpla}-\r{gcll} by the feedback interconnection of the plant
\eql{pla}{\arc
\mat{cc}{\t z\\y}=
\mas{c|ccc}{
\rho^{-1}I_d &\rho^{-1}I_d  &\rho^{-1}mI_d \hl
0            &-I_d          &L\Id-m\Id    \\
I_d          &0             &0             }
\mat{cc}{\t w\\u}}
with the controller
\eql{con}{
u=\mas{c|c}{\rho^{-1}\Aalp&\rho^{-1}\Balp\hl \Calp&\Dalp}y.}
Moreover,
the weighted interconnection \r{wcll} is clearly given by closing the loop with the same controller \r{con} and the following filtered version of \r{pla}:
\eql{wpla}{\arc
\mas{c|cc}{\Af&\Bf&0\hl \Cf(\La)&\Df(\La)&0\\0&0&I_d}
\mas{c|ccc}{
\rho^{-1}I_d &\rho^{-1}I_d  &\rho^{-1}mI_d \hl
0            &-I_d          &L\Id-m\Id    \\
I_d          &0             &0             }=
\mas{cc|ccc}{
\Af      &0            &-\Bf          &\Bf(L-m)       \\
0        &\rho^{-1}I_d &\rho^{-1}I_d  &\rho^{-1}mI_d  \hl
\Cf(\La) &0            &-\Df(\La)     &-\Df(\La)(L-m) \\
0        &I_d          &0             &0              }.
}
Recall that $\Af$ is Schur. Therefore, the structure of $\Acltf$ in \r{wcll} clarifies that
$\rho^{-1}(\Aal+\Bal m\Cal)$ is Schur iff  \r{con} stabilizes \r{wpla}.

For some given $\La\in\Lasr$, we can apply Theorem~\ref{Tsyn} in order to
characterize the existence of algorithm parameters that achieve \r{lmia} for some $\Xc$ as an LMI feasibility test.
However, the joint search over $\La\in\Lasr$ and $\red{X}$, $\Y$ in the resulting inequalities
remains non-convex.
This is a commonly encountered  problem in robust controller synthesis (see e.g.\ \cite{VeeSch14}), and has been also noted for algorithm design in the recent paper \cite{LesSei20}. 

Our progress over all existing results is to show how to perform such a {\bf simultaneous convex search} in \r{lmia} over the algorithm variables $\Aalp,\Balp,\Calp,\Dalp$, the multiplier parameters $\La$ and the Lyapunov matrix $\Xc$.

\remark{\label{remcon}
Once some $\La\in\Lasr$ has been determined, the corresponding algorithm parameters can be directly determined on the basis of Theorem~\ref{Tsyn} as in \cite{Gah94}. Let us include an important structural remark at this point. All matrices
involved in \r{wpla} and \r{QRS} do admit the Kronecker structure $M\otimes \Id$ with suitable matrices $M$ and
the standard Kronecker product ``$\otimes$''. This makes it possible to work w.l.o.g.\ with
$\red{X}$ and $\Y$ in Theorem~\ref{Tsyn} that admit such a structure, and the steps in \cite{Gah94}
generate algorithm parameters that inherit this structure as well. This so-called dimensionality reduction \cite{LesRec16}
implies that the computational complexity for algorithm design is independent of $d$, and that it generates
algorithms that can be applied for arbitrary dimesions $d\in\N$.
}

\subsection{Convexification of Operator Formulation}

Recall that the feedback interconnection of \r{pla} and \r{con} is called the closed-loop system. Let us start by showing that
all such closed-loops obtained by stabilizing controllers can be expressed as
\eql{Youla}{
\left\{T_1+T_2\Q T_3 \st \Q\in\cal\Q\right\}}
with the following set of stable systems:
\eql{Q}{
{\cal\Q}:=\left\{\Q=\mas{c|c}{\Aq&\Bq\hl\Cq&\Dq}\st \Aq\text{\ \ is Schur}\right\}.}
For this so-called Youla parametrization we rely on \cite{Fra87,ZhoDoy96} and emphasize that the results directly carry over from continuous-time to discrete-time systems.

\lemm{\label{Lpar}
The set of all systems \r{keysys} parameterized by $\Aalp,\Balp,\Calp,\Dalp$ 
and such that $\rho^{-1}(\Aal+\Bal m\Cal)$ is Schur is equal to \r{Youla} where $\si:=Lm^{-1}-1$,
\eqn{
T_1=
\mas{c|cc}{
0&-\rho^{-1} I_d\hl
\si I_d&-I_d},\ \ T_2=L\Id-m\Id,\ \
T_3=\mas{cc|c}{
0            &0    &-\rho^{-1}I_d \\
\rho^{-1}I_d &0    &\rho^{-1}I_d  \hl
0            &I_d  &0             }.
}
Moreover, the correspondence between the algorithm parameters and $\Q$ is given by
$$
\mat{c|c}{\Aalp&\Balp\hl \Calp&\Dalp}=
\mat{cc|c}{
\red{-}m\Dq     &\red{-}m \Cq   &\red{-}m\Dq       \\
\rho \Bq &\rho\Aq &\rho\Bq    \hl
\Dq      &\Cq     &\Dq-m^{-1}I_d }.
$$
}

\proof{In control we associate to a linear system \r{sysm} its so-called
transfer matrix $\Cal(\z I-\Aal)^{-1}\Bal+\Dal$, whose entries are real rational and proper functions in $\z$. A calculation shows that \r{pla} has the transfer matrix
$$\arc
P(\z)=\mat{cc}{P_{11}(\z)&P_{12}(\z)\\P_{21}(\z)&P_{22}(\z)}=
\mat{cc}{-\Id&L\Id-m\Id\\ \frac{1}{\rho \z-1}\Id&\frac{1}{\rho \z-1}m\Id}.
$$
With
$
M(\z)=\t M(\z):=\frac{\rho \z-1}{\rho \z}\Id,\ N(\z)=\t N(\z):=\frac{1}{\rho\z}m\Id
$
we then infer that $P_{22}(\z)=N(\z)M(\z)^{-1}$. With
$
X(\z)=\t X(\z):=\Id,\ \ Y(\z)=\t Y(\z):=-m^{-1}\Id
$
we get the so-called double B\'ezout identity
$$
\mat{cc}{\t X(\z)&-\t Y(\z)\\-\t N(\z)&\t M(\z)}\mat{cc}{M(\z)&Y(\z)\\N(\z)&X(\z)}=I_{2d}.
$$
This permits us to apply \cite[Theorem 1 in Section 4.5]{Fra87}.
Specifically, if defining
\eqn{
T_1(\z)=P_{11}(\z)+P_{12}(\z)M(\z)\t Y(\z)P_{21}(\z),\ \
\bar T_2(\z)=P_{12}(\z)M(\z),\ \ \bar T_3(\z)=\t M(\z)P_{21}(\z),
}
the set of all closed-loop transfer matrices that can be obtained with stabilizing
controllers for \r{pla} is given by
$T_1(\z)+\bar T_2(\z)Q(\z)\bar T_3(\z)$
where $Q(\z)$ varies in the set of all transfer matrices associated to the elements in $\cal \Q$.
Since $M(\z)$ is a multiple of the identity matrix, we infer
$\bar T_2(\z)Q(\z)=P_{12}(\z)Q(\z)M(\z)$ and hence
$T_1(\z)+\bar T_2(\z)Q(\z)\bar T_3(\z)=
T_1(\z)+T_2(\z)Q(\z)T_3(\z)$ with $T_2(\z):=P_{12}(\z)$, $T_3(\z)=M(\z)\t M(\z)P_{21}(\z)$. Now note that
$T_2(\z)=(L-m)\Id$,
\eql{TS}{
T_1(\z)=-\Id-\frac{1}{\rho \z}\si\Id\te{and}
T_3(\z)=\frac{\rho \z-1}{(\rho \z)^2}\Id,
}
which do indeed have the state-space representations as in the lemma.

According to \cite[Theorem 12.17]{ZhoDoy96}, the controller's transfer matrix which corresponds to
$Q(\z)$ is obtained by feedback of the plant with transfer matrix
$$\arc
\mat{cc}{Y(\z)X(\z)^{-1}&\t X(\z)^{-1}\\X(\z)^{-1}&\!\!\!\!{-}X(\z)^{-1}N(\z)}=\mat{cc}{-m^{-1}\Id&\Id\\\Id&\!\!\!\!{-}\frac{1}{\rho \z}m\Id}
$$
and $Q(\z)$. In case that $Q(\z)$ is the transfer matrix of an element in \r{Q},
an elementary calculation shows that the related controller has the state-space description
$$
\mas{cc|c}{
\red{-}\rho^{-1}m\Dq   &\red{-}\rho^{-1}m\Cq  &\red{-}\rho^{-1}m\Dq  \\
\Bq                    &\Aq                   &\Bq                   \hl
\Dq                    &\Cq                   &\Dq-m^{-1}I_d          }.
$$
Matching with \r{con} reveals the relation of $\Q$ with the algorithm parameters as claimed.  \epro
}

Next we note that $\Pi(\La)$, $T_1$, $T_2$, $T_3$ and $\Q$ are all stable, which implies the very same property for
\eql{wclly}{\Pi(\La)(T_1+T_2\Q T_3)}
due to \r{prod} and \r{sum}. Just by combining  Theorem \ref{Trs} with
Lemma \ref{Lpar}, we infer that there exists an algorithm which achieves $\rho$-convergence if there exist $\La\in\Lasr$ and $\Q\in\c\Q$ such that \r{wclly} is SNR.
The key step to convexity is the  parameter change $\Z:=\Pi(\La)T_2\Q$,
as shown in the next lemma.

\newcommand{\tLa}{\red{\tilde\La}}
\newcommand{\bLa}{\red{\bar\La}}

\lemm{\label{Lco}
Let $\bLa:=\mat{ccccc}{0&\cdots&0&I}$ and suppose
that $\bm{\La}$ with $\bLa\in \bm{\La}\subset\Las$ is convex.
Then the following statements are equivalent:
\enu{
\item
$\Pi(\La)(T_1+T_2\Q T_3)$ is SNR for some $\La\in\bm{\La}$, $\Q\in\c\Q$.
\item
$\Pi(\La)T_1+\Z T_3$ is SNR for some $\La\in\bm{\La}$, $\Z\in\c\Q$.}
}



\proof{If $\Pi(\La)(T_1+T_2\Q T_3)=\Pi(\La)T_1+\Pi(\La)T_2\Q T_3$ is SNR for some $\Q\in\c\Q$, it suffices to observe that
$\Z:=\Pi(\La)T_2\Q\in\c\Q$ since $\Z$ admits the description
$$
\mas{cc|c}{
\Af         &\Bf(L-m)\Cq       &\Bf(L-m)\Dq      \\
0           &\Aq          &\Bq         \hl
\Cf(\La)    &\Df(\La)(L-m)\Cq  &\Df(\La)(L-m)\Dq }
$$
and $\Af$ as well as $\Aq$ are Schur. Hence 1. implies 2.

To show that 2. implies 1., pick $\La\in\bm{\La}$, $\Z\in\c\Q$ such that
$\ab{w,(\Pi(\La)T_1+\Z T_3)w}_2\leq -\eps\|w\|_2^2$ for all $w\in l_2^d$ and some $\eps>0$.
In a first step, we slightly perturb $\La$ as $\La+\bar\eps\bLa$ in order render $\Df(\La+\bar\eps\bLa)=
\Df(\La)+\bar\eps\Df(\bLa)$ invertible; since $\Df(\bLa)=I$, this is indeed true for all sufficiently small $\bar\eps>0$.
Since $T_1$ is stable, its $l_2$-induced operator norm $\|T_1\|$ is finite.
By $\Pi(\bLa)=I$ we infer
for all $w\in l_2^d$ that
\eqn{
\ab{w,(\Pi(\La+\bar\eps\bLa)T_1+\Z T_3)w}_2=
\ab{w,\Pi(\La)T_1+\Z T_3)w}_2+
\bar\eps\ab{w,T_1w}_2
\leq(-\eps+\bar\eps \|T_1\|)\|w\|_2^2.
}
All this permits us to fix some small $\bar\eps\in(0,1)$ such that
$\tLa:=\La+\bar\eps\bLa\in\bm{\La}$ (convexity), $\Df(\tLa)$ is invertible and $\Pi(\tLa)T_1+\Z T_3$ stays SNR.
Therefore, $\Pi(\tLa)^{-1}$ exists and can be expressed as
$$\arc
\mas{c|c}{\t\Af&\t\Bf\hld \t\Cf&\t\Df}:=
\mas{c|c}{
\Af-\Bf\Df(\tLa)^{-1}\Cf&\Bf\Df(\tLa)^{-1}\hld
-\Df(\tLa)^{-1}\Cf&\Df(\tLa)^{-1}}.
$$
Next we show that $\t\Af$ is Schur. To this end we fix $\delta:=\frac{L-m}{2}\in(0,L-m)$ and the map $g(x):=\frac{1}{2}\delta x^Tx$.
Then Lemma \ref{Liqc} is valid for the full class $\Las$ \cite{MicSch20}, and thus as well for $\tLa$. If $\ch z\in l_2^d$, we note that \r{grag} just gives $\ch w=\delta \ch z$ and we get $\t z=(L-m)\ch z-\delta\ch z=\delta\ch z$ in \r{sit}; by Lemma \ref{Liqc} we hence conclude $\langle \Pi(\tLa)\ch z,\ch z\rangle_2\geq 0$ for all $\ch z\in l_2^d$; this shows
$\langle \Pi(\tLa)z,z\rangle=\ab{\Pi(\La)z,z}_2+\bar\eps\ab{\Pi(\bLa)z,z}_2\geq \bar\eps\|z\|_2^2$ for all $z\in l_2^d$. Since $\Af$
is Schur and $-\Pi(\La)$ with a state-space description in terms of $(\Af,\Bf,-\Cf(\tLa),-\Df(\tLa))$ is SNR,
Lemma \ref{Lpr} shows that there exists some $X\cg 0$ with
\eql{lmi1}{\arc
\bul^T
\mat{cc|cc}{X&0&0&0\\0&-X&0&0\hl 0&0&0&I\\0&0&I&0}
\mat{cc}{\Af&\Bf\\I&0\hld -\Cf(\tLa)& -\Df(\tLa)\\0&I}\cl 0.
}
Now we exploit again that $\Df(\tLa)$ is invertible and perform a congruence transformation of \r{lmi1} with
\eqn{
\mat{cc}{I&0\\-\Df(\tLa)^{-1}\Cf(\tLa)&\Df(\tLa)^{-1}}
\te{to get}
\bul^T
\mat{cc|cc}{X&0&0&0\\0&-X&0&0\hl 0&0&0&I\\0&0&I&0}
\mat{cc}{\t \Af&\t \Bf\\I&0\hl 0& -I\\\t\Cf&\t \Df}\cl 0.
}
By inspection, the left-upper block of this inequality reads $(\t\Af)^TX\t\Af-X\cl 0$.
Because of $X\cg 0$ we infer that $\t\Af$ is indeed a Schur matrix.

Since $T_2$ just is a real invertible matrix, we can define $\Q:=T_2^{-1}\Pi(\tLa)^{-1}\Z$.
We infer $\Q\in\c\Q$, again just by using \r{prod}. Moreover, $\Pi(\tLa)(T_1+T_2\Q T_3)=\Pi(\tLa)T_1+\Z T_3$ shows that
$\Pi(\tLa)(T_1+T_2\Q T_3)$ is SNR.
\epro }

Corollary~\ref{Can} combined with Lemma~\ref{Lpar}
and Lemma~\ref{Lco} for $\bm{\La}:=\Lasr$
leads to the following result.

\coro{With $\rho\in(0,1)$, there exists an algorithm whose convergence rate is bounded as
$\rwc\leq\rho$ if there exist
$\La\in\Lasr$ and $\Z\in\c\Q$ such that $\Pi(\La)T_1+\Z T_3$ is SNR.}

Both $\Lasr$ and $\c\Q$ are convex and $\Pi(\La)T_1+\Z T_3$ is affine in $\La$ and $\Z$. Since the SNR property is a convex constraint, we have shown that the algorithm design problem is indeed convex as a feasibility problem over the infinite dimensional space $\Lasr\times \c\Q$.

\subsection{Algorithm Synthesis by LMIs}\label{Salgsyn}

Testing whether there exist $\La\in\Lasr$ and $\Z\in\cal\Q$ for which $\Pi(\La)T_1+\Z T_3$ is SNR can even be turned into
a finite dimensional convex feasibility problem. Towards this end, we represent $\Z\in\c\Q$ as
\eql{con3}{
\Z=\mas{c|c}{\Az&\Bz\hl\Cz&\Dz}
}
and use \r{gpla}-\r{gcll} to see that
$\Pi(\La)T_1+\Z T_3$ results from the feedback interconnection of the plant \r{gpla}
with the controller $u=\Z y$ for
\equ{\label{sys3}
\mat{c|c|c}{A&B_1&B\hl C_1(\La)&D_1(\La)&E\hl C&F&0}:=
\arc
\mat{ccc|c|c}{
\Af           &\Bf\si       &0               &-\Bf                 &0    \\
0             &0            &0               &-\rho^{-1}\Id        &0    \\
0             &\Id          &0               &\rho^{-1}\Id         &0    \hl
\Cf(\La)      &\Df(\La)\si  &0               &-\Df(\La)            &L\Id-m\Id  \hl
0             &0            &\Id             &0                    &0     }.
}
This viewpoint permits us to derive an LMI solution for the algorithm synthesis problem based on
Theorem \ref{Tsyn}, our second main result. The relevant LMIs can be more compactly expressed by using
the selection matrix
$$
\sel:=\mat{cc}{I_{\dim(\Af)}&0\\0&I_d\\0&0}.
$$

\theo{\label{Tsynlmi}There exists some $\La\in\Lasr$ and a controller \r{con3} such that $\Az$ is Schur and
$\Pi(\La)T_1+\Z T_3$ is SNR iff there exist $\red{X}$ and $\La$ satisfying
\gat{\label{syn1}
\La\in\Lasr,\ \ \red{X}\cg 0\te{and}
\bul^T\mat{cc|cc}{\red{X}&0&0&0\\ 0&-\red{X}&0&0\hl  0&0&0&\Id\\0&0&\Id&0}
\mat{ccc}{A\sel &B_1\\\sel &0\hl C_1(\La)\sel&D_1(\La)\\0&\Id}\cl 0.
}
If the LMIs \r{syn1} are feasible, there exists an algorithm whose convergence rate is bounded as $\rwc\leq\rho$.
}

\proof{Let us abbreviate the interconnection \r{con3}-\r{sys3} as determined according to \r{gcll} by
\eql{cll3}{
z=\mas{c|c}{\Acltf&\Bcltf\hld \Ccltf(\La)&\Dcltf(\La)}w
}
(which is an abuse of notation since the matrices differ from those in \r{wcll}).
Since $B$ in \r{sys3} is zero, $\Az$ is Schur iff $\Acltf$ is Schur. Therefore, $\Az$ is Schur and \r{con3}-\r{sys3} is SNR iff there exists some $\Xc\cg 0$ with
\r{lmia}. To apply Theorem \ref{Tsyn}, we choose the basis matrices
$$\arc
U:=\mat{ccc}{I_{\dim(\Af)}&0&0\\0&I_d&0\\0&0&0\\0&0&I_d},\ \
V:=\mat{ccc}{I_{\dim(\Af)}&0&0\\0&I_d&0\\0&0&I_d\\0&0&0}
$$
of $\ker\mat{cc}{C&F}$, $\ker\mat{cc}{B^T&E^T}$, respectively.
Then the second LMI in \r{lmis1} is just identical to third one in \r{syn1}; since the first LMI in \r{lmis1} implies $\red{X}\cg 0$, ``only if'' follows directly.

To show ``if'', let \r{syn1} hold. By the particular choice of $V$, the third inequality in \r{lmis1}
simplifies to
$\red{Y}-A\red{Y}A^T\cg 0$.
Since $A$ is Schur, we can take $\red{Y_0}\cg 0$ with $\red{Y_0}-A\red{Y_0}A^T=I$ and thus
obtain for any $\al>0$ a solution $\al\red{Y_0}$ of the third LMI in \r{lmis1}. Since $\red{X}\cg 0$, we can certainly find some large $\al_0>0$ such that $\red{Y}=\al_0\red{Y_0}$ also satisfies the first LMI in \r{lmis1}. Applying theorem \ref{Tsyn} completes the proof.\epro}

Like for algorithm analysis, the left-hand side of \r{syn1} constitute LMI constraints on $\red{X}$ and $\La$.
Feasibility of these LMIs is equivalent to the existence of $\Aalp,\Balp,\Calp,\Dalp$, $\La\in\Lasr$ and $\Xc$ with \r{lmia}, which is the desired convexification of algorithm synthesis, one of the main goals of this paper.

Let us now establish that one can even eliminate the unknown $\red{X}$ in Theorem~\ref{Tsynlmi}.

\coro{\label{Cpick}
Let $R:=\sum_{k=0}^l \rho^k(\Lai{k}Lm^{-1})$
for fixed $\La\in\Lasr$. Then there exists $\red{X}$ with \r{syn1} iff
\eql{lmiLa}{\arc
\mat{cc}{
\frac{1}{1-\rho^2}(R^T+R)&R^T+\Lai{0}\\
R+\LaiT{0}&\LaiT{0}+\Lai{0}}\cg 0.
}
}
\proof{Note that the last LMI in \r{syn1} is a generalized Stein inequality
\eql{st1}{\arc
\bm{A}^T\red{X}\bm{A}-\bm{E}^T\red{X}\bm{E}+\bm{C}^TN\bm{C}\cl 0\text{\ where\ }
\mat{c}{\bm{A}\\\bm{E}\\\bm{C}}:=
\mat{ccccc}{
\Af           &\Bf\si              &-\Bf                  \\
0             &0                &-\frac{1}{\rho}\Id  \\
0             &\frac{1}{\rho}\Id              &\frac{1}{\rho}\Id               \hl
I_{\dim(\Af)}        &0                &0                     \\
0             &\Id              &0                     \\
0             &0                &0                     \hl
\Cf(\La)      &\Df(\La)\si         &-\Df(\La)             \\
0             &0                &\Id                   },\ N:=\arc\mat{cc}{0&\Id\\\Id&0}.
}
We start by determining a congruence transformation on \r{st1} in order to render
$\bm{A}$ diagonal. If $U$ is chosen to satisfy the Sylvester equation
\gat{\label{sy1}
\Af U-U(\mbox{$\frac{1}{\rho}$}\Id)+\Bf(1+\si)=0,
}
we indeed have
$$
P^{-1}\bm{A}Q={\bm{\bar A}}:=
\mat{ccc}{\Af&0&0\\0&\frac{1}{\rho} \Id&0\\0&0& -\Id},\
P^{-1}\bm{E}Q=\bm{E}\te{for}
P:=\mat{ccc}{I&U&{\Bf}\\0&\Id &\frac{1}{\rho}\Id \\0&0& -\frac{1}{\rho}\Id },\ 
Q:=\mat{ccc}{I&U&0\\0&\Id &0\\0&-\Id & \Id }.
$$
A congruence transformation of \r{st1} with $Q$ leads to
\eql{st2}{
\bm{\bar A}^T\red{\bar X}\bm{\bar A}-\bm{E}^T\red{\bar X}\bm{E}+\bm{\bar C}^TN\bm{\bar C}\cl 0
\te{where}\red{\t X}:=P^T\red{X}P,\ \
\bm{\bar C}:=\bm{C}Q=
\mat{ccccc}{
\Cf(\La)      &R       &-\Df(\La)             \\
0             &-\Id    &\Id                      },
}
and
$
R:=\Cf(\La)U+\Df(\La)(1+\si).
$
By \r{sy1} we have
$U=(\mbox{$\frac{1}{\rho}$}\Id-\Af)^{-1}\Bf(1+\si)$ and hence we get
$
R=[\Cf(\La)(\mbox{$\frac{1}{\rho}$}\Id -\Af)^{-1}\Bf+\Df(\La)](1+\si).
$
Due to \r{fil} and $\si+1=Lm^{-1}$, this indeed matches with the definition of $R$ in the corollary.

To show ``only if'' let $\red{X}$ satisfy \r{syn1}. Then
$\red{\t X}:=P^T\red{X}P\cg 0$ satisfies \r{st2}. Its
right-lower $2d\times 2d$-block is denoted as $\red{Y}$
and still positive definite. Canceling the first block row/column of \r{st2}
gives
\equ{\label{st3}\arc
\mat{ccccc}{
\frac{1}{\rho}I               &0                    \\
0                             &-I      }^T
\mat{cc}{\red{Y_{11}}&\red{Y_{12}}\\\red{Y_{21}}&\red{Y_{22}}}
\mat{ccccc}{
\frac{1}{\rho}I                  &0                    \\
0                                &-I      }
-
\mat{cc}{\red{Y_{11}}&0\\0&0}
+
\mat{cc}{
-R^T-R&R^T+\Df(\La)\\
R+\Df(\La)^T&-\Df(\La)^T-\Df(\La)}\cl0.}
With a sign-change in the off-diagonal blocks and if recalling $\Df(\La)=\Lai{0}$, this is equivalent to
\eql{st4}{\arc
\mat{cc}{
R^T+R-\frac{1-\rho^2}{\rho^2}\red{Y_{11}}
&R^T\!+\!\Lai{0}\!-\!\frac{1}{\rho}\red{Y_{12}}\\
R+\LaiT{0}-\frac{1}{\rho}\red{Y_{21}}&\LaiT{0}+\Lai{0}-\red{Y_{22}}}\cg 0.
}
If $H:=R^T+R-\frac{1-\rho^2}{\rho^2}\red{Y_{11}}\cg 0$ is the left-upper block in here, we infer from $\rho\in(0,1)$
that
$
\frac{1}{1-\rho^2}(R+R^T)-\frac{1}{\rho^2}\red{Y_{11}}=
\frac{1}{1-\rho^2}H\cg H.
$
By $\tiny\arc
\mat{cc}{
\frac{1}{\rho^2}\red{Y_{11}}&\frac{1}{\rho}\red{Y_{12}}\\
\frac{1}{\rho}\red{Y_{21}}&\red{Y_{22}}}\cg 0$, the inequality
\r{st4} hence implies \r{lmiLa}.

To prove the converse, let \r{lmiLa} hold and define
$$
\red{Y}:=
\mat{cc}{
\frac{\rho^2}{1-\rho^2}(R^T\!\!+\!R\!-\!\eps I)&\rho(R^T+\Lai{0})\\
\rho(R+\LaiT{0})&\LaiT{0}\!\!+\!\Lai{0}-\eps I};
$$
here we can choose so small $\eps>0$ that $\red{Y}$ is positive definite. Moreover,
$\red{Y}$ obviously satisfies \r{st4} and thus \r{st3}. Since $\Af$ is Schur, we can choose $\red{\bar X_0}\cg 0$ with $\Af^T\red{\bar X_0}\Af-\red{\bar X_0}=-I$. Let us then define
$
\red{\bar X}:=\diag(\al\red{\bar X_0},\red{Y})\cg 0\te{with}\al>0
$
and consider $\bm{\bar A}^T\red{\bar X}\bm{\bar A}-\bm{E}^T\red{\bar X}\bm{E}+\bm{\bar C}^TN\bm{\bar C}$;
its right-lower $2d\times 2d$-block equals \r{st3} and is, therefore, negative definite; moreover, $\al$ only affects the left-upper block of this matrix, which actually just equals
$\Af^T(\al\red{\bar X_0})\Af-(\al\red{\bar X_0})=-\al I$; therefore, we can fix a sufficiently large $\al>0$ to make sure that $\red{\bar X}\cg 0$ satisfies \r{st2}. Then $\red{X}:=P^{-T}\red{\bar X}P^{-1}\cg 0$ is a solution of \r{st1} and hence of \r{syn1}, which finishes the proof.\epro}

\remark{\label{rem3}
Note that Corollary \ref{Cpick} can be linked to a beautiful mathematical approach for
solving $H_\infty$- and SNR-synthesis problems based on Nevanlinna-Pick interpolation (see e.g. \cite{LimAnd88,Kim89}). In fact,
Theorem \ref{Tsynlmi} concerns $\Pi(\La)T_1+\Z T_3$ in which $T_3$ is a stable system with as many inputs as outputs.
The related SNR-synthesis problem is classically said to be of the one-block type. It is also known that the
unstable zeros of the transfer matrix of $T_3$ play a key role in characterizing its solvability. Due to
\r{TS}, these are given by
$
z_1=\rho^{-1}\te{and}z_2=\infty.
$
As it turns out after a simple computation, -\r{lmiLa} is nothing but the so-called Pick matrix
\eql{pick}{
\mat{cc}{
\frac{H(z_1)^*+H(z_1)}{1-\bar z_1^{-1}z_1^{-1}}&\frac{H(z_1)^*+H(z_2)}{1-\bar z_1^{-1}z_2^{-1}} \\[1.5ex]
\frac{H(z_2)^*+H(z_1)}{1-\bar z_2^{-1}z_1^{-1}}&\frac{H(z_2)^*+H(z_2)}{1-\bar z_2^{-1}z_2^{-1}} }
}
where $H$ denotes the transfer matrix of $\Pi(\La)T_1$ (with the definition in \cite{ByrGeo01} which permits zeros at infinity).
}
For the particular class of multipliers \r{Cpick}, we can even go one step further and explicitly characterize the set of those
parameters $\rho$ for which LMI \r{lmiLa} is feasible.

\coro{\label{Cexp}Let $\ka:=Lm^{-1}$ and $\le\geq 1$. Then there exists some $\La\in\Lasr$ with \r{lmiLa} iff $1-\frac{1}{\sqrt{\kappa}}<\rho$.}

\proof{In view of the Kronecker product structure of the elements in $\Lasr$ and homogeneity of \r{syn1}, we can fix  $\lai{0}=1$ and express \r{lmiLa} as
\eql{lmi2}{
\mat{cc}{
\frac{2\ka}{1-\rho^2}(1+\sum_{k=1}^l \rho^k\lai{k})&\bul\\
1+\ka(1+\sum_{k=1}^l \rho^k\lai{k})&2}\cg 0.
}
If we set $\al:=\frac{1+\rho}{1-\rho}>1$, we infer that
$
\mat{cc}{
\frac{2}{1-\rho^2}\be&1+\be\\
1+\be&2}\cg 0\te{iff}\be\in\left(\al^{-1},\al\right).
$
Moreover, since $\le\geq 1$, one can check that
$
\left\{
1+\sum_{k=1}^l \rho^k\lai{k}\st
1+\sum_{k=1}^l \rho^{-k}\lai{k}\geq 0,\ \lai{k}\leq 0\right\}
$
is the interval $[1-\rho^2,1]$; this set is compact and convex; the maximum is $1$ and the minimum is
computed by dualization of the corresponding linear program.
Taken together, \r{syn1} is feasible iff 
$
\left(\al^{-1},\al\right)\cap[\kappa(1-\rho^2),\kappa]
$
is not empty. The infimal $\rho\in(0,1)$ for which this is true is determined by the equation
$\frac{1+\rho}{1-\rho}=\kappa(1-\rho^2)$, which indeed gives $1-1/\sqrt{\kappa}$. \epro}

If using Zames-Falb multipliers of any length $\le\geq1$ to certify convergence,
Corollary~\ref{Cexp} means that $1-\frac{1}{\sqrt{\kappa}}$ is the optimal rate
that is achievable among all algorithms \r{nl}.
In view of \cite{ScoFre18}, this proves for the first time that the triple momentum algorithm
is guaranteed to be optimal even if allowing for Zames-Falb multipliers of length $\le>1$.
This also clarifies why various attempts
to improve the rate by manual tuning \cite{LesRec16} or sum-of-squares optimization \cite{FazMor18} of the algorithm parameters were not successful.
Our computation of an explicit optimal rate-bound for design is analogous to what has been achieved for the analysis of Nesterov's algorithm in \cite{SafJos18}. Our approach brings out the intrinsic system theoretic reasons for the limits of performance in algorithm design;
this holds for both the value of the optimal rate (determined by two zeros of some transfer matrix), and for the insight
that algorithms \r{nl} with matrices $\Aal$ of dimension larger than two are not beneficial. All this is a consequence
of systematically formulating algorithm design as a controller synthesis problem for the plant \r{pla}.

\subsection{An extension with a numerical example}\label{Sfmul}

We emphasize that our algorithm design approach is more powerful than just proving Corollary~\ref{Cexp}.
This is illustrated by following \cite{GraEbe20} and showing how one can exploit additional structural knowledge about the cost functions. Specifically,
for given matrices $M_f,L_f\in\R^{d\times d}$ with $0\cl M_f\cl L_f$, we consider the class
$\c F$ of functions $f\in C^2(\R^d,\R)$ satisfying
\eql{class}{M_f\cle \nabla^2f(x)\cle L_f\te{for all}x\in\R^d.}
One could take the triple momentum algorithm and achieve the convergence rate
\eql{ka}{1-\frac{1}{\sqrt{\ka}}\te{for}
\kappa=\frac{\la_{\max}(L_f)}{\la_{\min}(M_f)}.}
Instead, we can as well design algorithms based on the matrices $M_f$ and $L_f$ by solving a suitable LMI system.
For this purpose, we introduce
$
T:=(L_f-M_f)^{-\frac{1}{2}}
$
and define \r{sys3} by replacing
\eql{mLs}{
(m\Id,L\Id,\si)\te{with}(M_f,L_f,L_fM_f^{-1}-\Id).
}
We then arrive at the the following convex algorithm design result for the class $\c F$.
\theo{\label{Tsynlmi2}One can construct an algorithm which achieves the convergence rate
$\rho\in(0,1)$ for the class $\c F$ if there exist $\La\in T\Lasr T^T$ and
$\red{X}\cg 0$ which satisfy the LMI \r{syn1}.
}

\proof{We first observe that all insights in Section~\ref{Sstr} remain valid after the substitution \r{mLs}. We argue
that the same holds for Sections~\ref{Sana} and \ref{Siqc} with $\Gclt$ constructed based on \r{mLs}.
To this end let $f\in\c F$ be taken with $\nabla f(0)=0$ and
define
$g(z):=f(Tz)-\frac{1}{2}(Tz)^TM_f(Tz)$
to infer
$\nabla^2g(z)=T^T\nabla^2f(Tz)T-T^TM_fT$ and thus
$
0\cle \nabla^2 g(z)\cle  T^T(L_f-M_f)T=I\text{ for all }z\in\R^n.
$
Hence $g\in\c{S}_{0,1}$.
With the transformations $\bar  w:=T^Tw$ and $\bar  z:=T^{-1}z$
in \r{nl} we obtain
$$
\bar  w=\nabla g(\bar  z),\ \ \bar  z=\mas{c|c}{\Aal+\Bal M_f\Cal&\Bal T^{-T}\hld  T^{-1}\Cal&0}\bar  w.
$$
By just following the line of reasoning in Sections \ref{Sana} and \ref{Siqc}, Theorem~\ref{Trs} holds for the convergence rate $\rwc$ with respect to $\c F$ if replacing \r{keysys} with
$$
\mas{c|c}{\rho^{-1}(\Aal+\Bal M_f\Cal)&\rho^{-1}\Bal T^{-T}\hl (1-0)T^{-1}\Cal&-\Id};
$$
observe that we use $g\in\c{S}_{0,1}^0$ at this point. Since
$T^{-1}=T^T(L_f-M_f)$, this can be expressed as
$$
T^T\Gclt T^{-T}\te{with}\Gclt:=\mas{c|c}{\rho^{-1}(\Aal+\Bal M_f\Cal)&\rho^{-1}\Bal\hl (L_f-M_f)\Cal&-I_d}.
$$
Theorem~\ref{Trs} involves $\Pi(\La)T^T\Gclt T^{-1}$ being stable and SNR, which is equivalent to the same
conditions for the congruence transformed system $T\Pi(\La)T^T\Gclt$ and hence for $\Pi(T\La T)\Gclt$; here $T\Pi(\La)T^T=\Pi(T\La T^T)$ is shown with \r{fil} and a suitable state-coordinate change. All this reveals that Corollary \ref{Can} persists to hold for the class
$\c F$ if replacing $\Lasr$ with $T\Lasr T^T$. The proof is then concluded as that of Theorem \ref{Tsynlmi}.\epro
}
\begin{figure}\center
{\includegraphics[width=0.45\textwidth,trim=20 15 20 15,clip=true]{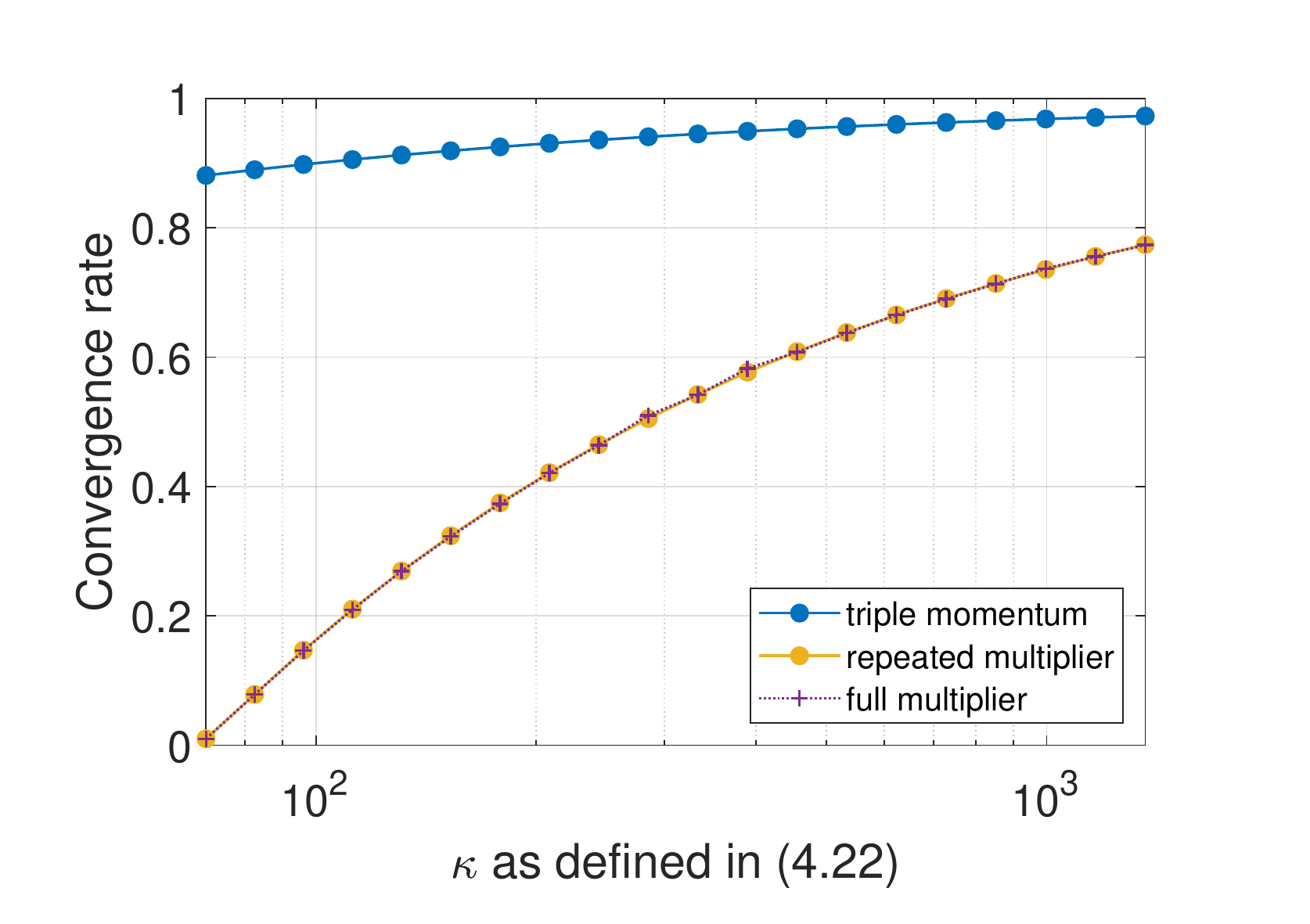}}
\caption{Optimal algorithm convergence rates versus condition number $\ka$ in \r{ka} for example in \cite{MicSch20}:
Triple momentum (blue), Theorem~\ref{Tsynlmi2} for $\Lasr$ (yellow) and Theorem~\ref{Tsynlmi2} for $\Las$ (purple).
}
\label{fign1}
\end{figure}
\tikzset{every path/.style={rounded corners}}

Once having determined some $\Laopt \in T\Lasr T^T$ for which the LMIs \r{syn1} in $\red{X}$ are feasible,
one can find related algorithm parameters $\Aalp,\Balp,\Calp,\Dalp$ as sketched in Section~\ref{Ssyna}. Moreover, the comments on dimensionality reduction carry over to the situation that $L_f-M_f$ (and hence $T$) are block-diagonal.

A concrete instance of the current setup are functions
$$f(x)=\frac{1}{2}x^TRx+h(Sx-s)\te{for}x\in\R^d$$ with given $R\in\R^{d\times d}$, $S\in \R^{e\times d}$, $s\in\R^e$ and any
$h\in\cS\cap C^2(\R^e,\R)$ where $R$ is positive definite.
Indeed, since $\nabla^2f(x)=R+S^T\nabla^2h(Sx-s)S$, we infer that \r{class} holds
with any small $\eps>0$ for
$$
M_f:=R+S^TmS-\eps I\te{and}L_f:=R+S^TLS.
$$
As motivated in \cite{MicSch20}, such cost functions  appear in model predictive control if handling
the constraint $Sx\leq s$ with a relaxed barrier function $b\in\cS\cap C^2(\R,\R)$ for the set $\{x\in\R\st x\leq 0\}$.
This results in the choice $h(x)=\sum_{i=1}^e b(x_i)$ with $\nabla h(x)=\col(b'(x_1),\ldots,b'(x_e))$, a so-called diagonally repeated nonlinearity. One can exploit this extra structure by using the full multiplier class $\Las$ instead of $\Lasr$
in Theorem~\ref{Tsynlmi2}; a proof relies on \cite[Theorem 9]{MicSch20} and the fact that Lemma~\ref{Lco} also applies to  $\bm{\La}=\Las$.
This offers yet another possibility for reducing conservatism in algorithm design.

We pick up the numerical example from \cite[Section 6.2]{MicSch20} for the latter class.
Figure~\ref{fign1} depicts the convergence rates of the triple momentum algorithm (blue), the structure exploiting algorithm from Theorem~\ref{Tsynlmi2} with repeated (yellow) and full multipliers (purple). In contrast to
\cite{MicSch20} (relying on non-convex design algorithms), we get identical rates for the two multiplier classes with our convex design algorithms.

\section{Generalization: Extremum Control}
\label{sec5}


\begin{figure}
\cen{
\begin{tikzpicture}[xscale=1,yscale=1]
\def\dl{5ex}
\node[sy3] (g1) at (0,0)  {$G_1$};
\coordinate[right=4*\dn of g1] (j);
\node[sy3,below=\dl of g1] (g2) at (0,0)  {$G_2$};

\node[sy3, left=\dl of g1] (f)  {$\nabla f$};
\node[sy3, left=3*\dl of g1] (k1)  {$K_1$};
\node[sy3, left=3*\dl of g2] (k2)  {$K_2$};
\node[sum, left= 3*\dn of k1] (s) {$+$};

\draw[->] (g1) -- node[]{$z$}  (f);
\draw[->] (f) --  node[]{$w$}  (k1);
\draw[->] (k1) --  (s);
\draw[->] (g2) -- node[pos=.84]{$v$} (k2);
\draw[->] (k2) -| (s);
\draw[-] (s.north) -- node[]{$u$} ([yshift=4*\dn] s.north) -| (j);
\draw[->] (j) -- (g1);
\draw[->] (j) |- (g2);

\end{tikzpicture}\hspace{10ex}
\begin{tikzpicture}[xscale=1,yscale=1]
\def\dl{5ex}

\node[sy3] (g) at (0,0)  {$G$};
\node[sy3,above =2*\dn of g] (f) {$\nabla f$};
\coordinate[right=4*\dn of g1] (j);
\node[sy3, left=1*\dl of g] (h)  {$H_d$};
\node[sy3, left=1*\dl of h] (k)  {$\bar K$};

\draw[->] (g) -- (h);
\draw[->] (h) -- node[]{$y$}  (k);
\draw[->] (k) -- node[]{$u$} ([xshift=-4*\dn] k.west) |-  node[pos=.95]{$z$} (f);
\draw[->] (f) -- node[pos=.5]{$w$} ([xshift=4*\dn] f.east) |- (g);


\end{tikzpicture}
\caption{Extremum Control: Optimization of the output of $G_1$ (left) or that of $\bar K$ (right).}
\label{fig2}
}
\end{figure}
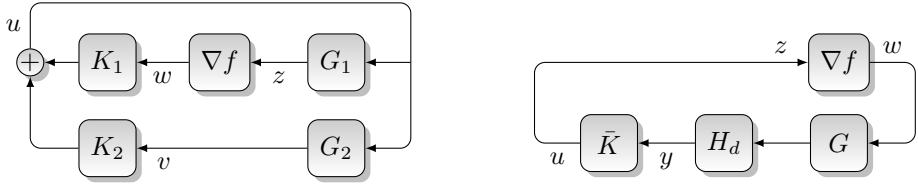

We now demonstrate that the proposed framework and the accompanying convexification result have a much wider scope than presented so far. They permit to systematically design optimization and learning algorithms with optimal convergence rates, even
with the presence of additional dynamics in the feedback loop.
Such dynamics may represent, for example, a model of a communication channels in optimization problems over networks, a noise filter if only a noisy gradient is available, the dynamics of a robot in a source seeking problem, or the dynamic properties of a hardware architecture like in neuromorphic computing.
A particularly nice scenario is extremum control as conceptually mentioned in \cite{AstWit89}. For a given system and a cost function, the goal is to design a controller that drives some system output to a steady-state in which the cost is minimal.


Among the many concrete instantiations of extremum control, we concentrate on the case where some linear system is given, the cost function $f$ is only known to belong to the class $\cS$, and the gradient of the cost function can be evaluated  \cite{MicEbe16a,NelMal18,LawSim18}.
To be concrete, we assume that the system is described as
\eql{syse}{
\mat{c}{z\\v}=\mat{c}{G_1\\G_2}u=\mas{c|c}{A_G&B_G\hl C_{G1}&D_{G1}\\C_{G2}&0}u,
}
with an input signal $u\in l^{n_u}$ used for control and two output signals
$z\in l^{d}$, $v\in l^{n_v}$ interpreted as follows. The first one is supposed to be asymptotically steered to $z_*\in\R^d$ with $\nabla f(z_*)=0$ for any cost $f\in\cS$.
The second output provides extra information about the system that can be exploited for control; it can be empty, which boils down to $n_v=0$.

The to-be-constructed dynamic controller (algorithm)  is a linear time-invariant system
that takes the two signals $w=\nabla f(z)$ and $v$ as its inputs and generates the control signal $u$ as its output:
\eql{cone}{
u=\mat{cc}{K_1&K_2}\mat{c}{w\\v}
,\te{}w=\nabla f(z).}
Altogether, \r{syse}-\r{cone} define the closed-loop system as depicted on the left in Fig.~\ref{fig2}.
With $\c G:=G_1(I-K_2G_2)^{-1}K_1$, a simple calculation shows that it can be described as
\eql{cG}{w=\nabla f(z),\ \ z=\c Gw.}
\remark{\l{rem4}Based on the state-space description of
\r{syse} and of $K_1$, $K_2$, one can calculate a representation \r{sysm} of $z=\c Gw$. We assume that
the ``$D$-matrix'' of  $K_1$ vanishes, which assures $\Dal=0$ in \r{sysm} and leads to \r{sys}.}

Given $\rho\in(0,1)$, the goal is to test whether there
exists a controller \r{cone} for \r{syse} which
achieves $\rho$-convergence for the interconnection \r{cG}
The infimum of all such $\rho$'s  is the extremum control rate $\rec$.

In Section \ref{sec_convexification} we have been only addressing the simple case $G_1=I$ and $G_2=[\,]$.
Despite the current more general setting, we are in the position to exploit the developed results in their full extent as seen next. The key is to adopt the generalized plant point-of-view.

\subsection{Setting up the Generalized Plant}

In the representation of $z=\c Gw$ by \r{sys} we can assume  w.l.o.g. that $(\Aal,\Cal)$ is detectable. By
Theorem \ref{Tnec}, we need to make sure that $\c G$ admits a factorization $\bar{\c G}H_d$
with $H_d$ being the integrator \r{int}. To enforce this structure we consider the feedback interconnection
\eql{gpe}{\arc
\mat{c}{z\hl y}=\mat{c|c}{
0&G_1 \hl
H_d&0   \\
0&G_2
}\mat{c}{w\hl u},\ \ u=\mat{cc}{\bar K_1&K_2}y
}
where the controllers $\bar K_1$ and $K_2$ can be freely chosen. 
With the abbreviation $\bar{\c G}:=G_1(I-K_2G_2)^{-1}\bar K_1$,
closing this loop indeed gives $z=\bar{\c G}H_dw$, and
\eql{constr}{K_1:=\bar K_1H_d}
assures the structure $\c G=G_1(I-K_2G_2)^{-1}K_1=\bar{\c G}H_d$.

We work again with state-space descriptions \r{gpla}-\r{gcon} of plant and controller in \r{gpe}, respectively.
Based on \r{syse}, the matrices for the plant can be taken as
\eql{gpp}{
\mat{ccc}{A&B_1&B\\ C_1&D_1&E\\ C&F&0}:=
\mat{cc|c|c}{
I &0      &I &0      \\
0 &A_G    &0 &B_G    \hl
0 &C_{G1} &0 &D_{G1} \hl
I &0      &0 &0      \\
0 &C_{G2} &0 &0      },
}
while those for the controller are free. Of course, then \r{gcll} leads to a state-space description of
\r{gpe} and, by inspection, $\Dcl$ indeed vanishes.
With the choice \r{constr}, this gives as well a state-space representation of $\c G$ in \r{cG}  which is
ensured to admit the structure \r{imp}.

For the purpose of synthesis, we need a plant-controller description for the corresponding system $\Gclt$ in \r{keysys}.
Recall that $\Gclt$ was obtained from \r{sys} by the signal transformations
\r{sit0}, \r{sit1} and \r{sit}.
If applied to the plant \r{gpla} with matrices \r{gpp}, these transformations lead to
\eql{gplat}{
\mat{c}{\t z\\ y}=
\mas{c|cc}{\At&\Bt_1&\Bt\hld \Ct_1&\Dt_1&\Et\\\Ct&\Ft&0}
\mat{c}{\t w\\u}
=
\mas{c|cc}{\rho^{-1}(A+B_1mC_1)&\rho^{-1}B_1&\rho^{-1}(B+B_1mE)\hl (L-m)C_1&-I&(L-m)E\\ C&0&0}
\mat{c}{\t w\\u}
}
corresponding to the former \r{pla}. A state-space description of $\Gclt$
is obtained from \r{gplat} interconnected with \r{con}.

Let us now formulate a test whether \r{gplat} admits a stabilizing controller which involves so-called uncontrollable and unobservable modes
\cite[Definition 3.6]{ZhoDoy96}.

\lemm{\label{Lstacon}
There exists a controller which stabilizes \r{gplat} iff
$(A_G,B_{G})$/$(A_G,C_{G})$ have no uncontrollable/unobservable modes in $\{\la\in\C\st |\la|\geq\rho\}$ and
$
\mat{cc}{A_G-I&B_{G}\\C_{G1}&D_{G1}}\text{\ has full row rank.\ }
$}
\proof{Since the triple $(\t A,\t B,\t C)$ is given by
$\arc
\left(
\rho^{-1}
\mat{cc}{
I &mC_{G1}     \\
0 &A_G   }
,
\rho^{-1}
\mat{cc}{
mD_{G1} \\
B_G     },
\mat{cc}{
I &0      \\
0 &C_{G2} }
\right),
$
one easily verifies with the so-called Hautus-test that the formulated conditions characterize that $(\t A,\t B)$/$(\t A,\t C)$ are stabilizable/detectable. This proves the claim by Remark \ref{rem2}.
\epro

This identifies necessary conditions on \r{syse} for achieving an extremum control rate of $\rho\in(0,1)$.
In particular, $G_1$ needs to be right invertible (requiring $d\leq n_u$) and should have no invariant zero at $1$ \cite{Won85a}.

\subsection{Convexification of Synthesis}

We assume that \r{gplat} admits a stabilizing controller.
Again, the key to convexification is the
description of all stabilized closed-loop systems as in \r{Youla}.
Now we follow a classical state-space procedure to construct this Youla-parameterization
\cite[Section 4.5]{Fra87}: Just  choose matrices $\Mt$ and $\Lt$ such that
$
\At+\Bt\Mt\text{ and }\At+\Lt\Ct\text{ are Schur}
$
and take
\equ{\arc\label{you}
\mat{cc}{T_1&T_2\\T_3&0}:=
\mas{cc|cc}{
\At +\Bt\Mt     &-\Bt\Mt    &\Bt_1        &\Bt \\
0               &\At+\Lt\Ct &\Bt_1+\Lt\Ft &0   \hld
\Ct_1+\Et\Mt    &-\Et\Mt    &\Dt_1        &\Et \\
0               &\Ct        &\Ft          &0   }.
}
Note that the left-hand side is a matrix of operators, whose blocks are defined by the one with a state-space representations on the right on the right-hand side.

Moreover, we exploit the structure of the multiplies $\Lasr$ in \r{Lasr} to arrive at the following result.

\lemm{\label{Lcoe}$\rho\in(0,1)$ satisfies $\rec\leq\rho$ if there exist $\La\in\Lasr$, $\Z\in\cal\Q$ s.th.
$\Pi(\La)T_1\!+\!T_2\Z T_3$ is SNR.}

\proof{
By \r{fil} and \r{Lasr}, the multiplier $\Pi(\La)$ admits the diagonal structure $\pi(\La)I_d$ with
$$
\pi(\La)=
\mas{cccc|c}{
0         &1           &\cdots        &0      &0        \\
\vdots    &\ddots      &\ddots        &\vdots &\vdots   \\
0         &\cdots      &0             &1      &0      \\
0         &\cdots      &0             &0      &1      \hl
\lai{\len}&\lai{\len-1}&\cdots        &\lai{1}&\lai{0} }
$$
having one input and one output only. Moreover, $T_2$ has the dimension $d\times n_u$ and
can be expressed as a operator matrix with entries $(T_2)_{ij}$ for $i=1,\ldots,d$, $j=1,\ldots,n_u$ that
are also systems with one input and one output. It is well-known that such systems commute as
$\pi(\La)(T_2)_{ij}=(T_2)_{ij}\pi(\La)$, which implies
$$\Pi(\La)T_2=(\pi(\La)\Id)T_2=T_2(\pi(\La)I_{n_u})=T_2\pi(\La).$$
Starting with
$\Pi(\La)T_1+T_2\Z T_3$ being SNR,
we follow the proof that 2. implies 1. in Lemma \ref{Lco}. Due to the structure of \ce{$\tLa$}, we infer $\Pi(\tLa)=\pi(\tLa)\Id$ and hence $\Pi(\tLa)^{-1}=\pi(\tLa)^{-1}\Id$. Therefore, $\Q:=\pi(\tLa)^{-1}\Z$ satisfies $\Q\in\Z$ and
assures that $\Pi(\tLa)(T_1+T_2\Q T_3)=\Pi(\tLa)T_1+T_2\Z T_3$ is SNR.
Again, the application of  Theorem~\ref{Trs} completes the proof.
\epro}

Verifying $\rec\leq\rho$ according to Lemma~\ref{Lcoe} is convex over $\Lasr\times\cal{\Q}$ and can be exactly turned into a finite dimensional LMI feasibility problem.
To see this we proceed as in Section \ref{Salgsyn} and express $\Pi(\La)T_1+T_2\Z T_2$ as
the interconnection
\eql{gplaew}{
\mat{c}{z\\y}=\mat{cc}{\Pi(\La)T_1&T_2\\T_3&0}\mat{c}{w\\u},\ \ u=\Z y.
}
Based on those of the multiplier \r{fil} and \r{you}, it is not difficult to construct a state-space representation of the plant in \r{gplaew}.
We dispense with writing down the matrices but note that these admit the structure
\eql{gps}{\arc
\mas{c|ccc}{A&B_1&B\hl C_1(\La)&D_1(\La)&E\\C&F&0}=
\mas{cc|cc}{
A_1          &0      &B_{11}      &0       \\
0            &A_2    &0           &B_{22}  \hl
C_{11}(\La)  &C_{12} &D_{1}(\La)  &E       \\
C_{21}       &0      &F           &0       }
}
in which $C_{11}(\La)$, $D_{1}(\La)$ and hence $C_1(\La)$ are affine in $\La$. Note that we abuse
notation since \r{gps} and \r{gpp} certainly are different plants.

We are now in the position to apply  (a specialized version of) a convexification procedure that has been first established in \cite{Sch00b}, see also \cite{RosSch20} for recent extensions.
The corresponding LMIs involve a symmetric decision variable $\red{W}$ with the same size and partition as $A$. Let us introduce the following functions in the variables $\red{W}$ and $\La$ (where we drop the arguments to save space):
\equ{\arc
\bm{W_1}\!:=\!\mat{cc}{\red{W_{11}}&\red{W_{12}}\\0&I},\ \bm{W_2}\!:=\!\mat{cc}{I&0\\-\red{W_{12}^T}&\red{W_{22}}},\label{var1}\ \
\red{\bm{Y}}:=\bm{W_1^T}\bm{W_2},\ \
\mas{c|ccc}{\bm{A}&\bm{B_1}\hl \bm{C_1}&\bm{D_1}}:=
\mas{c|ccc}{\bm{W_1^T}A\bm{W_2}&\bm{W_1^T}B_1\hl C_1(\La)\bm{W_2}&D_1(\La)}. 
}
It is crucial and easily checked by computation that the bold matrices all depend affinely  on $\red{W}$ and $\La$.

Moreover, let $U$, $\bm{V}$ be basis matrices of $\ker\mat{cc}{C&F}$, $\ker\mat{cc}{B^T&E^T}$ respectively; by \r{gps}
and with a basis $\col(V_1,V_2)$ of $\ker\mat{cc}{B_{22}^T&E^T}$, we can choose 
$
\bm{V}:=\mat{cc}{I&0\\0&V_1\\0&V_2}.
$


\theo{\label{Tsynlmi3}Either one of the following two equivalent conditions imply $\rec\leq\rho$ for $\rho\in(0,1)$:
\enu{
\item
There exists a $\La\in\Lasr$ and a controller \r{con3} such that $\Az$ is Schur and
$\Pi(\La)T_1+T_2\Z T_3$ is SNR.
\item There exist $\red{X}$, $\red{W}$ and $\La\in\Lasr$ satisfying $\mat{cc}{\red{X}&\bm{W_1}\\\bm{W_1^T}&\red{\bm{Y}}}\cg 0$ and
\eql{lmis1p}{\arc
\bul^T
\mat{cc|cc}{
\red{X}&0        &0 &0 \\
0      &-\red{X} &0 &0 \hl
0      &0        &0 &I \\
0      &0        &I &0 }
\mat{cc}{A&B_1\\I&0\hl C_1(\La)&D_1(\La)\\0&I}U\cl 0,\ \
\bul^T
\mat{cc|cc}{\red{\bm{Y}}&0&0&0\\0&-\red{\bm{Y}}^{-1}&0&0\hl 0&0&0&I\\0&0&I&0}
\mat{cc}{I&0\\-\bm{A^T}&-\bm{C_1^T}\hl 0&I\\ -\bm{B_1^T}&-\bm{D_1^T}}\bm{V}\cg 0.
}
}
}

%

{\em Sketch of proof.} By Lemma \ref{Lcoe} it suffices to show the equivalence of 1. and 2. To this end we characterize  1. with
Theorem \ref{Tsyn} for \r{QRS} in terms of a non-convex feasibility condition in the variables $\red{X}$, $\Y$ and $\La$. We then map $\Y\cg 0$ into
$
\mat{cc}{\red{W_{11}}&\red{W_{12}}\\\red{W_{21}}&\red{W_{22}}}:=
\mat{cc}{\red{Y_{11}^{-1}}&-\red{Y_{11}^{-1}}\red{Y_{12}}\\-\red{Y_{21}}\red{Y_{11}^{-1}}&\red{Y_{22}}-\red{Y_{21}}\red{Y_{11}^{-1}}\red{Y_{12}}}
$
with $\red{W_{11}}\cg 0$, $\red{W_{22}}\cg 0$. For \r{var1} one easily checks that
$\Y \bm{W_1}=\bm{W_2}$ and
$\det(\bm{W_1})\neq0$, $\det(\bm{W_2})\neq 0$.
Hence $\Y=\bm{W_2}\bm{W_1^{-1}}=\bm{W_1^{-T}}\bm{W_2^T}$. Thus,  by the third equation in \r{var1}, we get
\eql{subs}{
\Y=\bm{W_1^{-T}}\bm{\Y}\bm{W_1^{-1}}=\bm{W_2}\bm{\Y^{-1}}\bm{W_2^T}.}

Hence, the first inequalities  in \r{lmis1} and in 2. are related by a congruence transformation with the matrix $\diag(I,\bm{W_1})$.
Moreover, the second one in \r{lmis1} and the first in \r{lmis1p} are identical.

For the third inequality in \r{lmis1} we take the annihilator matrix
$V:=\diag(\bm{W_1},I)\bm{V}$.
This is fine since $B^T\bm{W_1}=B^T$ by the structure of $\bm{W_1}$, $B$, and hence
$\arc\mat{cc}{B^T&E^T}V=\mat{cc}{B^T\bm{W_1}&E^T}\bm{V}=\mat{cc}{B^T&E^T}\bm{V}=0$. An
inspection of the proof of Theorem \ref{Tsyn} in
\cite[Section 6.3]{Sch00siam} reveals that it causes no harm to take some $V$ which depends on $\Y$.
If we then substitute \r{subs} in  the third inequality of \r{lmis1}, we obtain
$$\arc
\bul^T
\mat{cc|cc}{\bm{\Y}&0&0&0\\0&-\bm{\red{Y^{-1}}}&0&0\hl 0&0&0&I\\0&0&I&0}
\mat{cc}{\bm{W_1^{-1}}&0\\-\bm{W_2^T}A^T&-\bm{W_2^T}C_1^T\hl 0&I\\ -B_1^T&-D_1^T}V\cg 0.
$$
With $V:=\diag(\bm{W_1},I)\bm{V}$ and \r{var1} we right away obtain the second inequality in \r{lmis1p}.
\epro

In fact, 2. constitutes convex constraints on all decisions variables; genuine LMIs are
obtained by taking a Schur complement w.r.t. $\red{\bm{Y}}^{-1}$ in the second inequality of \r{lmis1p}.

\subsection{Overall Design Procedure and Discussion}

Let us collect the steps to solve the extremum control problem for the system \r{syse}.
Note that this encompasses algorithm design for $G_1=\Id$ and $G_2=[\,]$ as in the first part of the paper.
\enu{
\item Choose some $\le\in\N_0$ in \r{fil} and $\rho\in(0,1)$. Set up the generalized plant \r{gplat} and
verify that the conditions in Lemma~\ref{Lstacon} hold.
\item Based on \r{gplat} and by using \r{you} and \r{fil}, construct a state-space representation of the generalized plant in
\r{gplaew} with the structure \r{gps}.
\item Check feasibility of the convex constraints in 2. of Theorem \ref{Tsynlmi3} (by solving the related LMI problem).
\item If feasible, the convergence rate $\rho$ is achievable for
a suitable controller \r{cone}  constructed as follows.
Choose the weighted version of the generalized plant \r{gplat} given as
\eql{wgplat}{\arc
\mas{cc|ccc}{
\Af         &\Bf\t C_1         &\Bf\t D_1         &\Bf\t E         \\
0           &\t A              &\t B_1            &\t B            \hld
\Cf(\Laopt) &\Df(\Laopt)\t C_1 &\Df(\Laopt)\t D_1 &\Df(\Laopt)\t E \\
0           &\t C              &\t F              &0               }.
}
For \r{wgplat} and \r{QRS}, the LMIs in Theorem \ref{Tsyn} are feasible.
With a related controller \r{gcon} and in view of \r{con}, \r{constr}, we infer that
$\arc
\mat{cc}{K_1&K_2}=
\mas{c|c}{\rho\Aco&\rho\Bco\hl \Cco&\Dco}\mat{cc}{H_d&0\\0&I_{n_v}}
$
achieves $\rho$-convergence for \r{syse}-\r{cone}.
}
}
By bisection over $\rho$, one can determine the infimal worst-case convergence rate $\rho_*(\le)$ that is achievable with some algorithm. Note that this best rate depends on the length $\le$ of the Zames-Falb multiplier \r{fil}.

We emphasize that any papers in the literature revolve around the case $\le=0$ leading to a multiplier without dynamics, which is related to the so-called circle-criterion or a version of the small-gain theorem. For example in \cite{MicEbe16a}, controllers are assumed to have an observer structure incorporating an integrator, and the design is split up into sequential observer and state-feedback synthesis steps relying on the small-gain theorem, both of which are typically conservative. The paper \cite{NelMal18} fixes a control structure and is confined to stability analysis for $\le=0$ only.
Dynamic multipliers are generally known to be considerably more powerful (see e.g. \cite{VeeSch16a} for analysis and \cite{Sch15} for synthesis). The test for $\le=1$ termed off-by-one circle criterion attracted special attention in \cite{LesRec16,LesSei20,ScoFre18} for algorithm analysis (see also Corollary~\ref{Cexp}).

Our approach overcomes various of these limitations in that we perform direct output-feedback synthesis and employ dedicated dynamic multipliers for the general extremum control problem.
Our numerical examples illustrate that it can be beneficial to work with multipliers of length larger than $1$, and that one can analyze the achievable convergence rates depending on suitable system theoretic properties of \r{syse}.

We have not made any assumptions on $G_1$ so far. In case that $G_1$ admits the structure $G_1=g_1I_{n_u}$ with a single input single output system $g_1$, we emphasize that the extensions as described in Section  \ref{Sfmul} go through with ease in the current more general setting. If, in addition, $G_2$ is empty and we work with multipliers that admit a Kronecker structure, the possibility for dimensionality reduction carries over as well (see Remark \ref{remcon}). The synthesis procedure will lead to a controller $K_1$ that also admits this structure. Concretely, this captures the optimal synthesis of an algorithm for the minimization of $f$ where
the information sent to the gradient first needs to pass a communication channel that is modeled by $g_1$.

\begin{figure}\center
\includegraphics[width=0.45\textwidth,trim=20 15 20 15,clip=true]{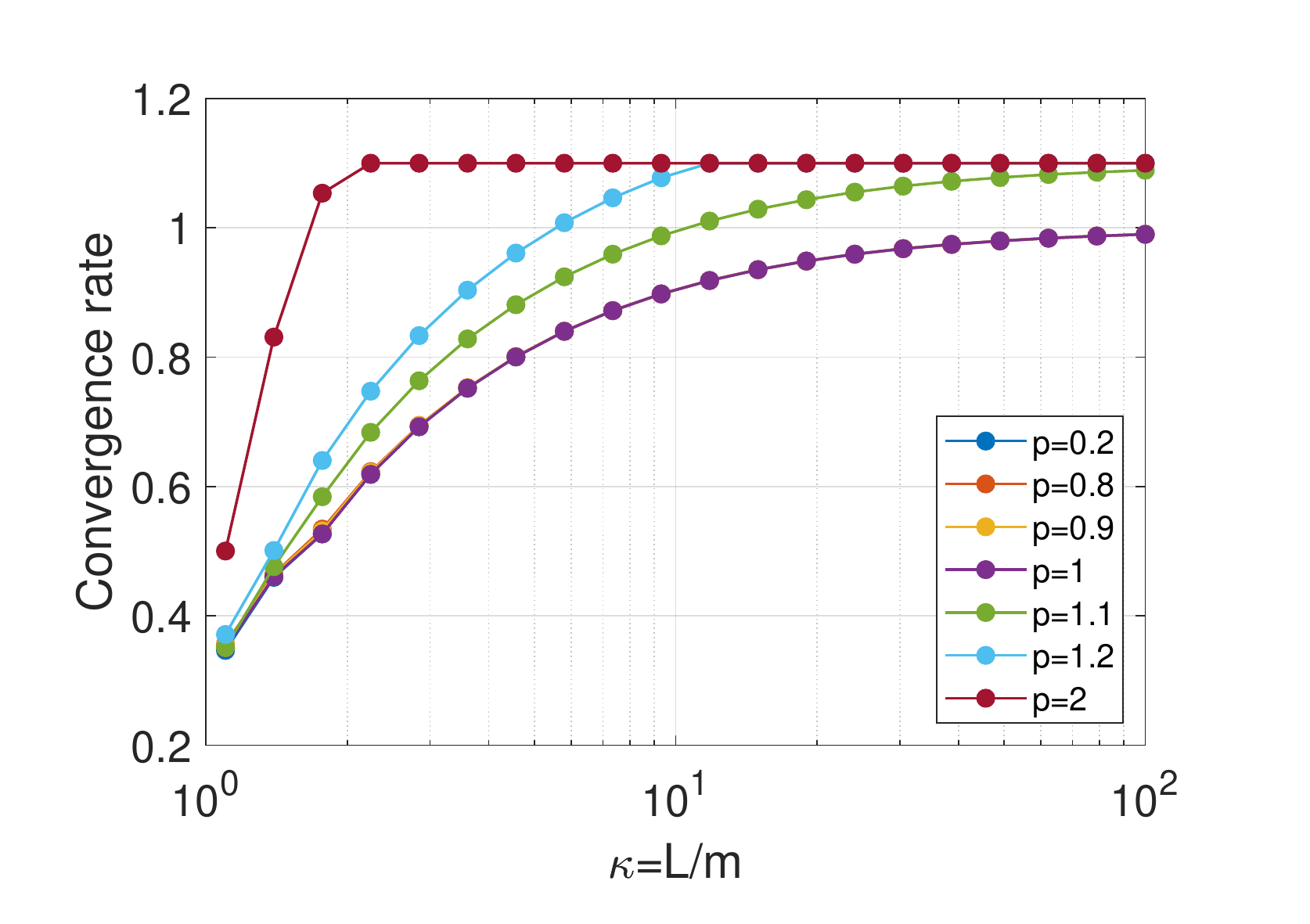}
\includegraphics[width=0.45\textwidth,trim=20 15 20 15,clip=true]{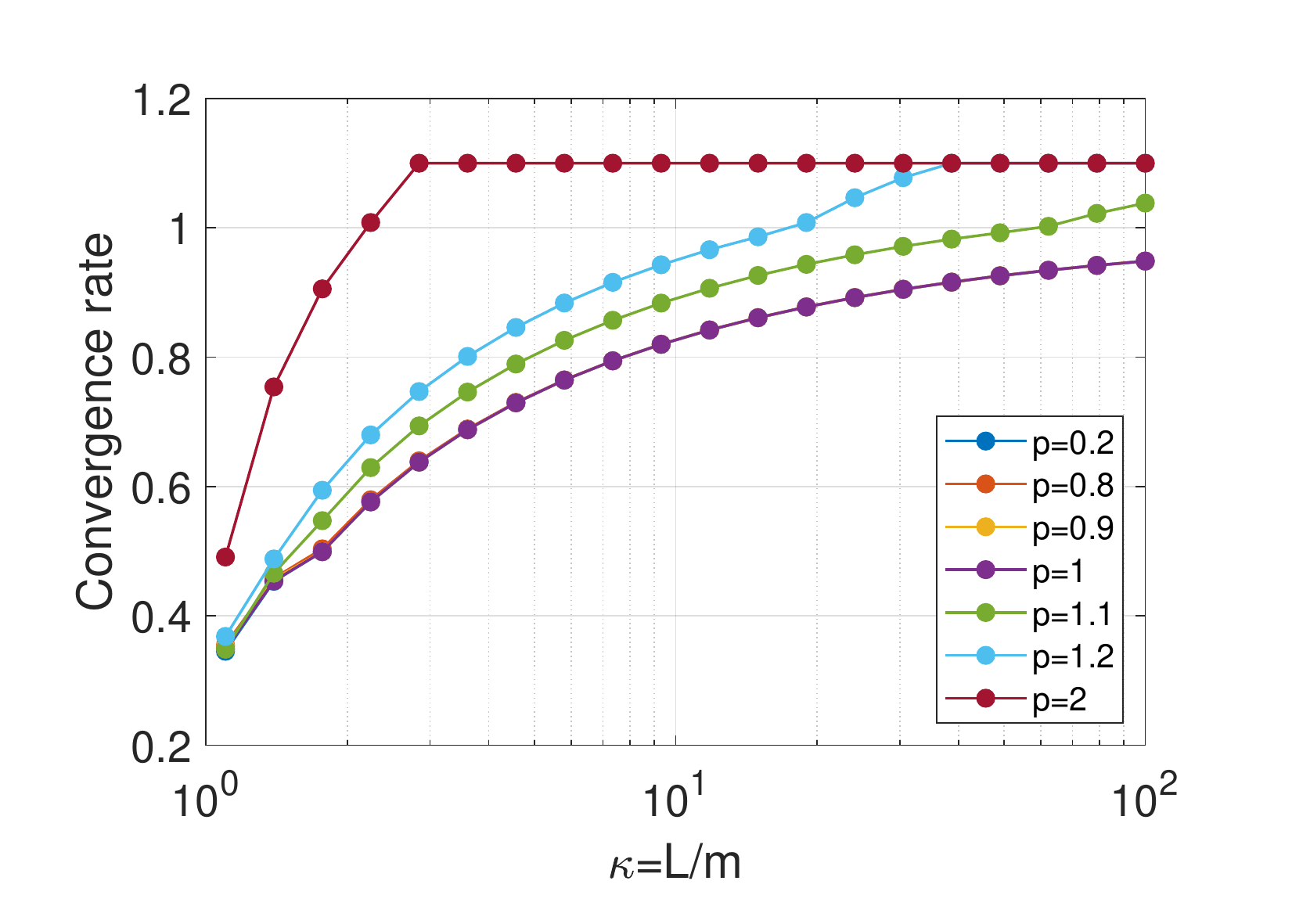}
\caption{Guaranteed rates for the example in Section \ref{Sex1} with
 $l=0$ (left) and $l=2$ (right).
}
\label{figex1}
\end{figure}

%

\subsection{Numerical examples}\label{Sexa}

\subsubsection{Example 1}\label{Sex1}
Relating the limits of performance of some controlled systems to properties of the underlying uncontrolled one is a classical research topic in control \cite{Zam81}. Our tools put us in the position to explore such limits of performance expressed by the achievable optimal convergence rate.
For a numerical illustration, we choose a very simple configuration with
$G_2=[\,]$ and a family of systems $G_1$ that admit the transfer functions
$
G_1(z)=\frac{z-0.5}{(z+p)(z+0.5)}
$
with the pole $p$ varying in $\{0.2, 0.8, 0.9, 1, 1.1, 1.2, 2\}$. The optimal achievable rates for static  ($l=0$) and dynamic multipliers ($l=2$) depending on the condition number $L/m$ are plotted  in Fig.~\ref{figex1}; note that the saturation at $1.1$ is due to the initialization of the bisection over $\rho$ with the interval $[0,1.1]$. The curves indicate a uniform improvement of the achievable rates if stepping from $l=0$ to $l=2$, but they do not improve any more for $l>2$. It is as well interesting to observe that the rates do not change when moving the pole $p$ inside the unit disk towards the boundary, but that they do get worse if $p$
moves further into the unstable region. Instead of discussing other interesting aspects of such trade-offs and fundamental performance limitations for the example, we conclude by emphasizing that the key aspect is the mere ability
to generate such plots also for many other scenarios.

%
%
%
%
%
\subsubsection{Example 2}\label{Sex2}
\begin{figure}\center
{\includegraphics[width=0.45\textwidth,trim=20 15 20 15,clip=true]{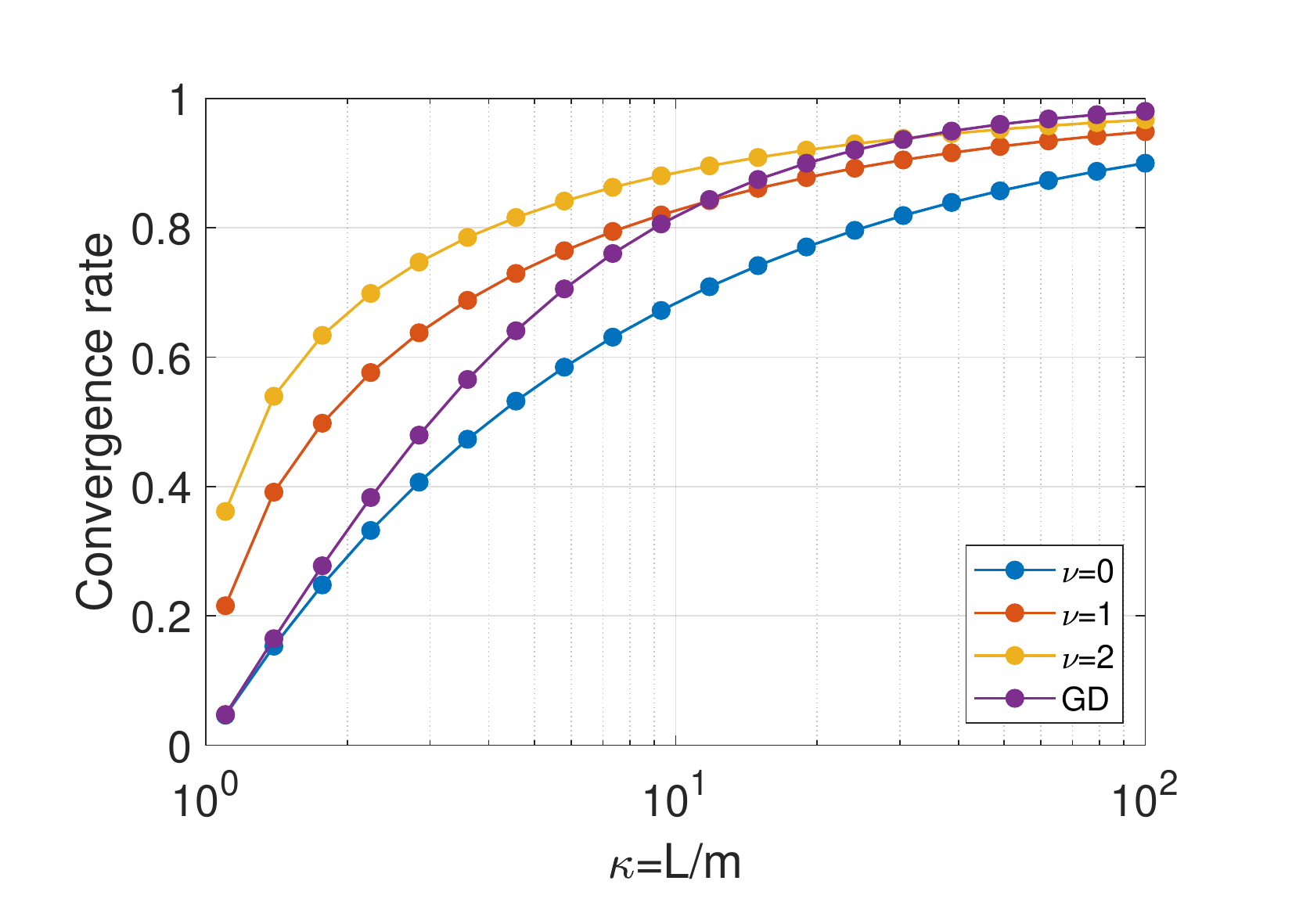}}
\includegraphics[width=0.45\textwidth,trim=20 15 20 15,clip=true]{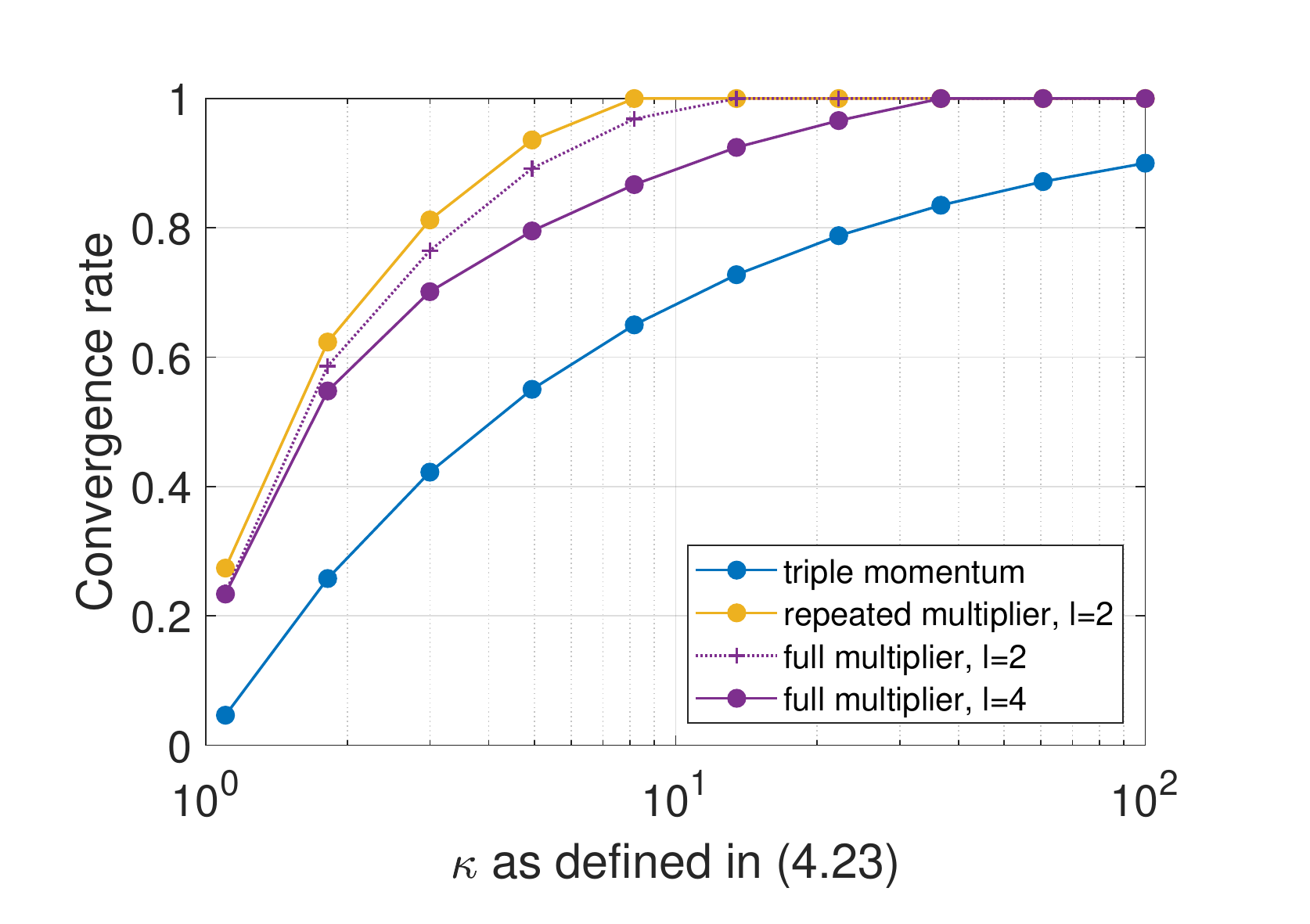}
\caption{Guaranteed rates for the example in Section \ref{Sex2} (left) and Section \ref{Sex3} (right).}
\label{fig3r}
\end{figure}

A particularly  interesting case is optimization with delayed gradient information,
as it appears in parallel optimization or optimization over networks, see e.g. \cite{ZheMen16}.
Hereby, convergence rates of gradient descent algorithms described as $x_{k+1} = x_k - \al \nabla f(x_{k-\nu})$ with step-size $\al>0$
and lag $\nu \in \mathbb{N}$ are studied. Such a delay in general accelerated gradient descent algorithms can be easily captured
by taking $G_1(z)=\frac{1}{z^\nu}$ and $G_2=[\,]$ on the left in  Fig. \ref{fig2}.
The achievable guaranteed convergence rates for $\nu=0,1,2$ are shown on the left in Fig. \ref{fig3r}.
Notice that $\nu=0$ corresponds to the triple-momentum method. As expected, a longer lag
leads to a lower performance. It is interesting to observe that our approach
allows to design an accelerated algorithm with delayed gradients that outperforms the standard gradient descent algorithm (GD) without delay for larger values of $\kappa$.

\subsubsection{Example 3}\label{Sex3}
Finally, we are not tied to the left configuration in Fig.~\ref{fig2}.
For example, one might require to optimize the rate of convergence
to an optimal steady-state $u_*$ of the controller's output on the right in Fig.~\ref{fig2} (in which the integrator $H_d$ is already displayed explicitly). Let us also assume that
$\nabla f$ is diagonally repeated as $\nabla f(z)=\col(b'(z),\ldots,b'(z))$ with any $b\in\cS\cap C^2(\R,\R)$ (see Sec.~\ref{Sfmul}).
Then algorithm synthesis can be convexified along the discussed lines for both repeated and full dynamic multipliers
$\Pi(\La)$ with $\La\in\Lasr$ and $\La\in\Las$, respectively.
If we pick
$\footnotesize\arc
G(z)=\mat{cc}{1&\frac{z-0.2}{z+1.2}\\\frac{z-1.2}{z+0.5}&1},
$
Fig.~\ref{fig3r} on the right depicts the achievable guaranteed rates for
repeated multipliers of length $l=2$ and for full multipliers of lengths $l=2$ and $l=4$.
These results not only reveal the benefit of
higher order dynamics in the multipliers, but also that of exploiting
the structure in the cost function.

\section{Conclusions}

Expanding on the point of view discussed in \cite{WanEli11}, it is one of our key messages that the generalized plant
view-point adopted in this paper offers otherwise unachievable conceptual and
structural insights into the analysis and synthesis of optimization algorithms.
We believe that this view-point has a high potential to stimulate further research at the interface of systems theory, optimization and machine learning. For example, it seems highly promising to incorporated recent advances in structured controller synthesis (see \cite{RosSch20} and references therein) for the convex design of distributed optimization algorithms.
Finally, we believe that our approach can be utilized in algorithm design problems
where performance properties such as the mitigation of the sensitivity against noise come into play.



\begin{thebibliography}{10}

\bibitem{AndVon73}
B.~Anderson and S.~Vongpanitlerd.
\newblock {\em {N}etwork {A}nalysis and {S}ynthesis}.
\newblock Prentice Hall, Englewood Cliffs, New Jersey, 1973.

\bibitem{AstWit89}
K.~{\AA}str{\"o}m and B.~Wittenmark.
\newblock {\em Adaptive Control}.
\newblock Addison-Wesley series in electrical engineering : control
  engineering. Addison-Wesley, 1995.

\bibitem{ByrGeo01}
C.~Byrnes, T.~Georgiou, and A.~Lindquist.
\newblock {A} generalized entropy critertion for {N}evanlinna-{P}ick
  interpolation with degree constraint.
\newblock {\em IEEE Trans. Aut. Control}, 46(5):822--839, 2001.

\bibitem{DesVid75}
C.~Desoer and M.~Vidyasagar.
\newblock {\em {F}eedback Systems: {I}nput-Output Approach}.
\newblock Academic Press, London, 1975.

\bibitem{DueEbe12}
H.-B. D\"urr and C.~Ebenbauer.
\newblock On a class of smooth optimization algorithms with applications in
  control.
\newblock {\em {IFAC} Proceedings Volumes}, 45(17):291--298, 2012.

\bibitem{FazMor18}
M.~Fazlyab, M.~Morari, and V.~M. Preciado.
\newblock Design of first-order optimization algorithms via sum-of-squares
  programming.
\newblock In {\em 2018 {IEEE} Conference on Decision and Control ({CDC})},
  2018.

\bibitem{Fra87}
B.~Francis.
\newblock {\em {A} course in {$H_\infty$} control theory}.
\newblock Springer-Verlag, Berlin, 1987.

\bibitem{Gah94}
P.~Gahinet.
\newblock {A} new parametrization of {$H_\infty$} suboptimal controllers.
\newblock In {\em International Journal of Control}, pages 1031--1051, 1994.

\bibitem{GahApk94}
P.~Gahinet and P.~Apkarian.
\newblock {A} linear matrix inequality approach to {$H_\infty$} {C}ontrol.
\newblock {\em Int. J. Robust Nonlin.}, 4:421--448, 1994.

\bibitem{GraEbe20}
D.~Gramlich, C.~Ebenbauer, and C.~W. Scherer.
\newblock Convex synthesis of accelerated gradient algorithms for optimization
  and saddle point problems using {L}yapunov functions.
\newblock {\em arXiv:2006.09946 [math.OC]}.

\bibitem{HorJoh85}
R.~Horn and C.~Johnson.
\newblock {\em {M}atrix {A}nalysis}.
\newblock Cambrigde University Press, New York, 1985.

\bibitem{HuLes17}
B.~Hu and L.~Lessard.
\newblock Control interpretations for first-order optimization methods.
\newblock In {\em 2017 {American} {Control} {Conference} ({ACC})}, pages
  3114--3119, May 2017.
\newblock ISSN: 2378-5861.

\bibitem{IwaSke94}
T.~Iwasaki and R.~Skelton.
\newblock {All} controllers for the general {${\cal H}_\infty$} control
  problem: {LMI} existence conditions and state space formulas.
\newblock {\em Automatica}, 30:1307--1317, 1994.

\bibitem{Kim89}
H.~Kimura.
\newblock Conjugation, interpolation and model-matching in ${H}_\infty$.
\newblock {\em International Journal of Control}, 49(1):269--307, 1989.

\bibitem{LawSim18}
L.~S.~P. Lawrence, J.~W. Simpson-Porco, and E.~Mallada.
\newblock Linear-convex optimal steady-state control.

\bibitem{LesRec16}
L.~Lessard, B.~Recht, and A.~Packard.
\newblock {A}nalysis and {D}esign of {O}ptimization {A}lgorithms via {I}ntegral
  {Q}uadratic {C}onstraints.
\newblock {\em SIAM Journal on Optimization}, 26(1):57--95, 2016.

\bibitem{LesSei20}
L.~Lessard and P.~Seiler.
\newblock Direct synthesis of iterative algorithms with bounds on achievable
  worst-case convergence rate.
\newblock In {\em 2020 American Control Conference ({ACC})}. {IEEE}, jul 2020.

\bibitem{LimAnd88}
D.~J.~N. Limebeer and B.~D.~O. Anderson.
\newblock An interpolation theory approach to {H}8 controller degree bounds.
\newblock {\em Linear Algebra and its Applications}, 98:347--386, Jan. 1988.

\bibitem{MasOha98}
I.~Masubuchi, A.~Ohara, and N.~Suda.
\newblock {LMI}-based controller synthesis: a unified formulation and solution.
\newblock {\em Int. J. Robust Nonlin.}, 8:669--686, 1998.

\bibitem{MegRan97}
A.~Megretski and A.~Rantzer.
\newblock {S}ystem analysis via {I}ntegral {Q}uadratic {C}onstraints.
\newblock {\em IEEE T. Automat. Contr.}, 42:819--830, 1997.

\bibitem{MicEbe16a}
S.~Michalowsky and C.~Ebenbauer.
\newblock Extremum control of linear systems based on output feedback.
\newblock In {\em 2016 {IEEE} 55th Conference on Decision and Control ({CDC})}.
  {IEEE}, 2016.

\bibitem{MicSch20}
S.~Michalowsky, C.~Scherer, and C.~Ebenbauer.
\newblock Robust and structure exploiting optimisation algorithms: an integral
  quadratic constraint approach.
\newblock {\em International Journal of Control}, pages 1--24, 2020.

\bibitem{NelMal18}
Z.~E. Nelson and E.~Mallada.
\newblock An integral quadratic constraint framework for real-time steady-state
  optimization of linear time-invariant systems.
\newblock In {\em 2018 Annual American Control Conference ({ACC})}, 2018.

\bibitem{Nes18}
Y.~Nesterov.
\newblock {\em Lectures on Convex Optimization}, volume 137 of {\em Springer
  Optimization and Its Applications}.
\newblock Springer International Publishing, 2018.

\bibitem{Pol87}
B.~Polyak.
\newblock {\em Introduction to Optimization}.
\newblock Optimization Software, Inc., New York, 1987.

\bibitem{Ran96}
A.~Rantzer.
\newblock On the {K}alman-{Y}akubovich-{P}pov lemma.
\newblock {\em Systems {\&} Control Letters}, 28(1):7--10, 1996.

\bibitem{RosSch20}
C.~A. Rosinger and C.~W. Scherer.
\newblock A flexible synthesis framework of structured controllers for
  networked systems.
\newblock {\em {IEEE} Transactions on Control of Network Systems}, 7(1):6--18,
  2020.

\bibitem{SafJos18}
S.~Safavi, B.~Joshi, G.~Franca, and J.~Bento.
\newblock An explicit convergence rate for nesterov's method from {SDP}.
\newblock In {\em 2018 {IEEE} International Symposium on Information Theory
  ({ISIT})}. {IEEE}, jun 2018.

\bibitem{Sch15}
C.~Scherer.
\newblock {G}ain-scheduling control with dynamic multipliers by convex
  optimization.
\newblock {\em SIAM J. Contr. Optim.}, 53(3):1224--1249, 2015.

\bibitem{SchGah97}
C.~Scherer, P.~Gahinet, and M.~Chilali.
\newblock Multiobjective output-feedback control via {LMI} optimization.
\newblock {\em {IEEE} Transactions on Automatic Control}, 42(7):896--911, 1997.

\bibitem{Sch00b}
C.~W. Scherer.
\newblock {D}esign of {S}tructured {C}ontrollers with {A}pplications.
\newblock In {\em Proc. 39th IEEE Conf. Decision and Control}, Sydney,
  Australia, 2000.

\bibitem{Sch00siam}
C.~W. Scherer.
\newblock {R}obust {M}ixed {C}ontrol and {LPV} {C}ontrol with {F}ull {B}lock
  {S}calings.
\newblock In L.~{El~Ghaoui} and S.~Niculescu, editors, {\em Advances in Linear
  Matrix Inequality Methods in Control}, pages 187--207. SIAM, Philadelphia,
  2000.

\bibitem{SchWei99}
C.~W. Scherer and S.~Weiland.
\newblock {\em {L}inear matrix inequalities in control}.
\newblock Lecture Notes, Delft University of Technology, 1999.

\bibitem{ScoFre18}
B.~V. Scoy, R.~A. Freeman, and K.~M. Lynch.
\newblock The fastest known globally convergent first-order method for
  minimizing strongly convex functions.
\newblock {\em IEEE Control Systems Letters}, 2(1):49--54, 2018.

\bibitem{TayDro21}
A.~Taylor and Y.~Drori.
\newblock An optimal gradient method for smooth (possibly strongly) convex
  minimization.
\newblock {\em arXiv:2101.09741 [math.OC]}.

\bibitem{VeeSch14}
J.~Veenman and C.~Scherer.
\newblock {A} synthesis framework for robust gain-scheduling controllers.
\newblock {\em Automatica}, 50(11):2799--2812, 2014.

\bibitem{VeeSch16a}
J.~Veenman, C.~W. Scherer, and H.~K\"{o}ro\u{g}lu.
\newblock {R}obust stability and performance analysis based on integral
  quadratic constraints.
\newblock {\em European Journal of Control}, 31:1--32, 2016.

\bibitem{WanEli11}
J.~Wang and N.~Elia.
\newblock A control perspective for centralized and distributed convex
  optimization.
\newblock In {\em IEEE Conference on Decision and Control and European Control
  Conference Orlando, FL, USA}. IEEE, 2011.

\bibitem{Won85a}
W.~Wonham.
\newblock {\em {L}inear {M}ultivariable {C}ontrol}.
\newblock Springer-Verlag, Berlin, 3rd edition, 1985.

\bibitem{Zam81}
G.~Zames.
\newblock Feedback and optimal sensitivity: Model reference transformations,
  multiplicative seminorms, and approximate inverses.
\newblock {\em {IEEE} Transactions on Automatic Control}, 26(2):301--320, 1981.

\bibitem{ZheMen16}
S.~Zheng, Q.~Meng, T.~Wang, W.~Chen, N.~Yu, Z.-M. Ma, and T.-Y. Liu.
\newblock Asynchronous stochastic gradient descent with delay compensation.

\bibitem{ZhoDoy96}
K.~Zhou, J.~Doyle, and K.~Glover.
\newblock {\em {R}obust and {O}ptimal {C}ontrol}.
\newblock Prentice Hall, Upper Saddle River, New Jersey, 1996.

\end{thebibliography}
\end{document}